\documentclass[fontsize=12pt,a4paper,headings=normal,
twoside=false,leqno,parskip=half-,abstract=true]{scrartcl}
\usepackage[english]{babel}
\usepackage[utf8]{inputenc}
\usepackage{hyperref}
\hypersetup{
 pdftitle={Bifurcation from infinity: a compactification approach},
 pdfauthor={Phillipo Lappicy},
 colorlinks=true,
 linkcolor=blue,
 citecolor=blue,
 filecolor=blue,
 urlcolor=blue}

\usepackage{graphicx}
\usepackage{pgfplots}
\usepackage[format=plain,labelfont=bf,font=small]{caption}
\usepackage{xcolor}
\usepackage[arrow, matrix, curve]{xy}
\usepackage{float}

\usepackage{caption}
\usepackage{comment}
\usepackage[textwidth=2cm,textsize=small,backgroundcolor=none]{todonotes}

\usepackage{tabulary}
\usepackage{array}

\usepackage{amsmath,amsthm}
\usepackage{amssymb} 
\usepackage{latexsym}

\usepackage{tikz}

\newtheorem{thm}{Theorem}[section]
\newtheorem{lem}[thm]{Lemma}

\newenvironment{pf}[1][Proof]{\begin{trivlist}
\item[\hskip \labelsep {\bfseries #1}]}{\end{trivlist}}

\usepackage[notref,notcite,color,final 
]{showkeys}

\begin{document}

\title{{\huge{On bifurcation from infinity: \\ a compactification approach
}}}

\author{
 \\
{~}\\
José M. Arrieta*, Juliana Fernandes$^\circ$ and Phillipo Lappicy*\\ 
}

\date{ }
\maketitle
\thispagestyle{empty}

\begin{abstract}
We consider a scalar parabolic partial differential equation on the interval with nonlinear boundary conditions that are asymptotically sublinear.
As the parameter crosses critical values (e.g. the Steklov eigenvalues), it is known that there are large equilibria that arise through a bifurcation from infinity (i.e., such equilibria converge, after rescaling, to the Steklov eigenfunctions). 
We provide a compactification approach 
to the study of such unbounded bifurcation curves of equilibria, their stability, and heteroclinic orbits. 
In particular, we construct an induced semiflow at infinity such that the Steklov eigenfunctions are equilibria. 
Moreover, we prove the existence of infinite-time blow-up solutions that converge, after rescaling, to certain eigenfunctions that are equilibria of the induced semiflow at infinity.
%
%

\ 


\ 

\end{abstract}


\vfill

$\ast$\\
{\footnotesize
Dept. Análisis Matemático y Matemática Aplicada, 
Universidad Complutense de Madrid, Spain\\
Instituto de Ciencias Matemáticas (ICMAT), CSIC-UAM-UC3M-UCM, Spain\\
\texttt{arrieta@mat.ucm.es} and \texttt{philemos@ucm.es}\\
} 


{\footnotesize
$\circ$\\
Instituto de Matemática, Universidade Federal do Rio de Janeiro, Brazil\\
\texttt{jfernandes@im.ufrj.br}\\
}

\newpage
\pagestyle{plain}
\pagenumbering{arabic}
\setcounter{page}{1}

\section{Introduction}\label{sec:intro}

\numberwithin{equation}{section}
\numberwithin{figure}{section}
\numberwithin{table}{section}

Nonlinear boundary conditions play an important role in modelling real life phenomena and can be used as control terms to acquire certain (in)stability properties of solutions. 
Hence, the interplay between the dynamics induced by the boundary conditions and by the interior kinetics is of fundamental importance. The goal of this paper is to provide a comparison of the effects of boundary nonlinearities.

We consider the semilinear parabolic equations
\begin{subequations}\label{intro:PDE}
\begin{align}
    u_t&=u_{xx}-u, \qquad\qquad\qquad\quad\,\,\,\, x\in (0,1), \quad t>0,\\   
    -u_x(t,0)&=\lambda u(t,0) + g(u(t,0)), \qquad\qquad\qquad\qquad t>0,\\
    u_x(t,1)&=\lambda u(t,1) + g(u(t,1)), \qquad\qquad\qquad\qquad t>0,
\end{align}
\end{subequations}
with initial data $u(0,x)=u_0(x)$ in a suitable phase-space $X$. 
We assume $g(u)$ is Lipschitz, odd, monotone, bounded, $g(0)=0$ and $g'(0)=1$; an example is $g(u)=\arctan(u)$. 

Under these hypotheses, for any initial data $u_0(x)\in X:=L^2([0,1])$, there is a unique solution $u(t,\cdot )\in C^0([0, \infty),X)$ of equation~\eqref{intro:PDE}, which is global and depends continuously on the initial data and on the parameter $\lambda\in \mathbb{R}$; see \cite{ACRB99,AC00,ACRB00} and \cite[Section 7]{ARRB07}. 
Indeed, equation~\eqref{intro:PDE} can be rewritten in the following abstract formulation:
\begin{subequations}\label{intro:PDE2}
\begin{align}
    u_t&=u_{xx}-u+G (u), \qquad\qquad\quad x\in (0,1), \,\quad t>0,\\   
    -u_x(t,0)&=\lambda u(t,0), \qquad\qquad\qquad\qquad\qquad\qquad\qquad t>0,\label{intro:linBC0}\\
    u_x(t,1)&=\lambda u(t,1), \qquad\qquad\qquad\qquad\qquad\qquad\qquad t>0,\label{intro:linBC1}
\end{align}
\end{subequations}
with initial condition $u(0,x)=u_0(x)\in X$, for an appropriate $G(u)$ that is also bounded, see \cite{ACRB99,AC00,ACRB00}. 
Considering the operator $A=\partial_x^2 - Id$ with linear Robin conditions \eqref{intro:linBC0}-\eqref{intro:linBC1}, this amounts to that solutions of~\eqref{intro:PDE2} satisfy the variation of constants formula: 
\begin{equation}\label{VoC}
    u(t,x)= \mathrm{e}^{A t}u_0(x) + \int_0^t \mathrm{e}^{A (t-s)} G(u(s,x)) ds.
\end{equation}
In fact, we will show that the global semiflow generated by equation~\eqref{intro:PDE} possess a global attractor $\mathcal{A}_\lambda$ with different features depending on $\lambda\in\mathbb{R}$; see Lemma~\ref{lem:nondiss}. On one hand, the semiflow may be dissipative and possess a compact attractor; on the other hand, the semiflow may be non-dissipative with an unbounded attractor. 
Either way, the attractor is the union of unstable manifolds of all equilibria (i.e. time-independent solutions): 
\begin{equation}\label{char:A}
    \mathcal{A}_\lambda = \bigcup_{u_*\in\mathcal{E}^b} W^u (u_*),
\end{equation}
where $\mathcal{E}^b$ is the set of all bounded equilibria of \eqref{intro:PDE}. Such characterization of the global attractor through the unstable manifold of the equilibria can be seen in \cite{Henry81} for the dissipative case and in \cite{FernandesBortolan} for the non-dissipative case.
This follows from the gradient structure of the semiflow, due to the Lyapunov function
\begin{equation}
    E(u)=\frac{1}{2}\int_0^1 |u_x|^2+|u|^2 - \lambda \left( \frac{u^2(t,x)}{2} + \int_0^{u(t,x)} g(s)ds\right)\bigg|_{x =0}^{x =1},
\end{equation}
and thus the $\alpha$- and $\omega$-limit sets of bounded solution of \eqref{intro:PDE} 
consist of equilibria, see~\cite{LBeatriz}.

{
\color{black}
Our goal is to describe $\mathcal{A}_\lambda$ explicitly.  
It is known that there are suitably large equilibria of~\eqref{intro:PDE}, which are generated by a resonant mechanism at the boundary, see \cite{ARRB07,ARRB09}. 
Indeed, there are bifurcation curves from infinity, 
see \cite{Rabinowitz,Toland,Stuart,Dancer,Toland2} for pioneering work on these matters.
We will provide a compactification approach to analyze such unbounded behavior of the semiflow~\eqref{VoC}, akin to \cite{Hell09,BenGal10,RochaPimentel15,LP}. 
In particular, we will show that the unbounded bifurcation curves arise from \emph{equilibria at infinity}, since the semiflow~\eqref{VoC} induces a projected dynamics at infinity consisting of equilibria and heteroclinics. Moreover, we also obtain the existence of infinite-time blow-up solutions, which can be regarded as \emph{heteroclinics to infinity}. 
}

The remaining of the paper is organized as follows. In Section~\ref{sec:ubdd}, we explicitly construct the unbounded bifurcation curves of equilibria and analyze their properties. 
In Section~\ref{sec:att}, we describe bounded and unbounded (infinite-time blow-up) solutions of~\eqref{intro:PDE}, and thereby whenever there is a compact or unbounded attractor. 
In Section~\ref{sec:P}, we provide the Poincaré projection for equation~\eqref{intro:PDE}. In particular, we show that the unbounded bifurcating branches of equilibria converge to an equilibria of the projected semiflow at infinity; we also discuss the existence of heteroclinics at infinity. In Section~\ref{sec:Thm}, we sum up our main results. In Section~\ref{sec:disc}, we discuss our approach and propose open problems.

\section{Main results}

\subsection{Unbounded curves of equilibria}\label{sec:ubdd}

There are bounded equilibria of \eqref{intro:PDE}, if $\lambda$ is not a solution of the problem:
\begin{subequations}\label{intro:Steklov}
\begin{align}
    0&=\Phi_{xx}-\Phi, \qquad\qquad x\in (0,1),\label{intro:Steklov1}\\    
    -\Phi_x(0)&=\sigma \Phi(0),\label{intro:Steklov2}\\
    \Phi_x(1)&=\sigma \Phi(1),\label{intro:Steklov3}
\end{align}
\end{subequations}
which is called the Steklov eigenvalue problem, \textcolor{black}{see~\cite{Stekloff}}.
To find the Steklov eigenvalues and eigenfunctions, note that general solutions of \eqref{intro:Steklov1} are given by $\Phi(x)=A \mathrm{e}^x + \tilde{A}\mathrm{e}^{-x} $, where we can find $A,\tilde{A}\in\mathbb{R}$ through the boundary conditions \eqref{intro:Steklov2}-\eqref{intro:Steklov3}: 
\begin{subequations}\label{abSteklov}
\begin{align}
    -A+\tilde{A}&=\sigma (A+\tilde{A}),\\
    A\mathrm{e} - \tilde{A}\mathrm{e}^{-1}&=\sigma (A\mathrm{e} + \tilde{A}\mathrm{e}^{-1}).
\end{align}
\end{subequations}
The first equation implies that $\tilde{A}=A(1+\sigma)/(1-\sigma)$, and plugging this into the second equation amounts to $(\sigma+1)^2-\mathrm{e}^2(\sigma-1)^2=0$.
        %
%
Therefore, there are only two Steklov eigenvalues, which solve the linear problem~\eqref{intro:Steklov}, given by 
%
%
\begin{equation}
    \sigma_1:=\frac{\mathrm{e}-1}{\mathrm{e}+1} \approx 0.462\ldots \qquad \text{ and } \qquad \sigma_2:=\frac{1}{\sigma_1}=\frac{\mathrm{e}+1}{\mathrm{e}-1}\approx 2.163\ldots.
\end{equation}
Hence, each Steklov eigenvalue possess an associated Steklov eigenfunction such that $\tilde{A}_1=A_1(1+\sigma_1)/(1-\sigma_1)=A_1\mathrm{e}$ and $\tilde{A}_2=A_2(1+\sigma_2)/(1-\sigma_2)=-A_2\mathrm{e}$, i.e. 
\begin{equation}
    \Phi_1(x):=A_1[\mathrm{e}^x+\mathrm{e}^{1-x}], \qquad \Phi_2(x):=A_2[\mathrm{e}^x-\mathrm{e}^{1-x}],
\end{equation}
where $A_1,A_2$ can be chosen in order to normalize the eigenfunctions in any desirable norm. For example, $A_1:=1/(1+\mathrm{e}),A_2:=1/(1-\mathrm{e})$ normalize $\Phi_1(x),\Phi_2(x)$ in $L^\infty$.
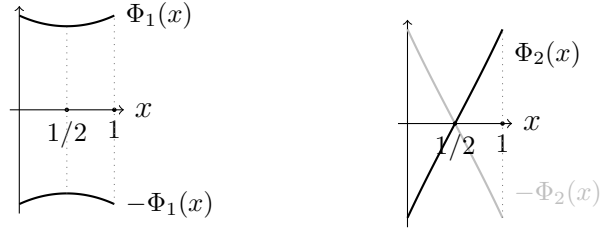
\begin{figure}[H]\centering
    \begin{tikzpicture}[scale=1.25]
    \draw[->] (-0.1,0) -- (1.1,0) node[right] {$x$};
    \draw[->] (0,-1.1) -- (0,1.1);

    \draw [thick, domain=0:1,variable=\t,smooth] plot ({(\t)},{(exp(\t)+exp(1-\t))/(1+exp(1))})  node[anchor=west]{\footnotesize{$\Phi_1(x)$}}; 
    \draw [thick, domain=0:1,variable=\t,smooth] plot ({(\t)},{-(exp(\t)+exp(1-\t))/(1+exp(1))})  node[anchor=west]{\footnotesize{$-\Phi_1(x)$}}; 


    \filldraw [black] (0.5,0) circle (0.5pt) node[anchor=north]{\footnotesize{$1/2$}};
    \filldraw [black] (1,0) circle (0.5pt) node[anchor=north]{\footnotesize{$1$}};
    
    \draw[gray,dotted] (1,-1) -- (1,1);
    \draw[gray,dotted] (0.5,-0.886) -- (0.5,0.886);

    \end{tikzpicture}
    \hspace{2cm}
    \begin{tikzpicture}[scale=1.25]
    \draw[->] (-0.1,0) -- (1.1,0) node[right] {$x$};
    \draw[->] (0,-1.1) -- (0,1.1);


    \draw [lightgray,thick, domain=0:1,variable=\t,smooth] plot ({(\t)},{(exp(1-\t)-exp(\t))/(exp(1)-1)}) node[anchor=south west]{\footnotesize{$-\Phi_2(x)$}}; 
    \draw [thick, domain=0:1,variable=\t,smooth] plot ({(\t)},{-(exp(1-\t)-exp(\t))/(exp(1)-1)}) node[anchor=north west]{\footnotesize{$\Phi_2(x)$}}; 

    \filldraw [black] (0.5,0) circle (0.5pt) node[anchor=north]{\footnotesize{$1/2$}};
    \filldraw [black] (1,0) circle (0.5pt) node[anchor=north]{\footnotesize{$1$}};
    
    \draw[gray,dotted] (1,-1) -- (1,1);

    \end{tikzpicture}
\caption{The two Steklov eigenfunctions $\Phi_1,\Phi_2$, which possess the symmetries $\Phi_1(x)=\Phi_1(1-x)$ and $\Phi_2(x)=-\Phi_2(1-x)$. 
        In particular, the minimum of $\Phi_1(x)$ and the zero of $\Phi_2(x)$ occur at $x=1/2$. 
}\label{Fig:EF}
\end{figure}
Each Steklov eigenvalue amounts to a bifurcation curve from infinity, possessing a positive and a negative branch, which can be computed explicitly. For the existence of branches see~\cite[Theorem 3.4 and 4.5]{ARRB07}, where the subcritical or supercritical character is discussed in \cite[Theorem 3.4 and 3.5]{ARRB09}, their stability in \cite[Proposition 7.1 and 7.3]{ARRB09}, and their continuation 
in \cite[Theorem 3.3]{ARRB07}. \textcolor{black}{See also~\cite{Cushing,StuartToland}}.
%
%

We expect that solutions in each bifurcation branch has the same symmetry as the respective Steklov eigenfunction, see Figure~\ref{Fig:EF}. The branch of solutions that emanates from $\Phi_1(x)$, denoted by $u^1(x)$, should be invariant under $u^1(x)\mapsto u^1(1-x)$; and the branch of solutions that emanates from $\Phi_2(x)$, denoted by $u^2(x)$, should satisfy $u^2(x)\mapsto -u^2(1-x)$.

On one hand, we look for equilibria of \eqref{intro:PDE} 
given by $u^1(x)=c_1 [\mathrm{e}^x + \mathrm{e}^{1-x}] $, where $c_1\in\mathbb{R}$ can be found by plugging such solution to the boundary conditions, yielding 
\begin{equation}\label{c_1}
    -(1-\mathrm{e})c_1=\lambda (1+\mathrm{e})c_1 + g((1+\mathrm{e})c_1).
\end{equation}
Hence, we can solve this explicitly
\begin{equation}
    \lambda(c_1)=\sigma_1 - \frac{g((1+\mathrm{e})c_1)}{(1+\mathrm{e})c_1}.
\end{equation}
%
%
On another hand, we look for equilibria of \eqref{intro:PDE} given by $u^2(x)=c_2 [\mathrm{e}^x - \mathrm{e}^{1-x}] $, where $c_2\in\mathbb{R}$ can be found through the boundary conditions, which for $g$ odd yields: 
\begin{equation}\label{c_2}
    (\mathrm{e}+1)c_2=\lambda (\mathrm{e}-1)c_2 + g((\mathrm{e}-1)c_2).
\end{equation}
Again, this can be solved explicitly
\begin{equation}
    \lambda(c_2)=\sigma_2 - \frac{g((\mathrm{e}-1)c_2)}{(\mathrm{e}-1)c_2}.
\end{equation}
%
%
Thus, by construction, we obtain that each of the branches converge to the corresponding normalized Steklov eigenfunction, i.e., for each $i=1,2$, we have that
\begin{equation}\label{EQconvSteklov}
    \frac{u_i(x)}{||u_i(x)||_{L^\infty(\partial (0,1))}} \xrightarrow{\lambda\to \sigma_i} \pm \Phi_i (x), \qquad \text{ in } C^\beta([0,1]), \beta\in (0,1).
\end{equation}
This was known in \cite{ARRB07,ARRB09}, which is in contrast with the approximation scheme in \cite{BCP}.

\begin{figure}[H]\centering
    \begin{tikzpicture}[scale=0.75]
    \node (1) at (1, 3) {\boxed{\text{\footnotesize{$g(u):=\arctan(u)$}}}};
    
    \draw[->] (-0.75,0) -- (2.5,0) node[right] {\footnotesize{$\lambda$}};
    \draw[->] (0,-2) -- (0,2);

    \draw [thick, domain=-2:2,variable=\t,smooth,rotate=-90] plot ({(\t)},{(e-1)/(e+1)-(rad(atan((e+1)*\t)))/((e+1)*\t)}) node[above left] {\footnotesize{$-c_1(\lambda)$}};
    \draw [thick, domain=-2:2,variable=\t,smooth,rotate=-90] plot ({(\t)},{(e+1)/(e-1)-(rad(atan((e-1)*\t)))/((e-1)*\t)}) node[above right] {\footnotesize{$-c_2(\lambda)$}};
    \draw [thick, domain=2:-2,variable=\t,smooth,rotate=-90] plot ({(\t)},{(e-1)/(e+1)-(rad(atan((e+1)*\t)))/((e+1)*\t)}) node[below left] {\footnotesize{$+c_1(\lambda)$}};
    \draw [thick, domain=2:-2,variable=\t,smooth,rotate=-90] plot ({(\t)},{(e+1)/(e-1)-(rad(atan((e-1)*\t)))/((e-1)*\t)}) node[below right] {\footnotesize{$+c_2(\lambda)$}};

    \filldraw [black] (0.462-1,0) circle (2pt) node[anchor=north]{\footnotesize{$\sigma^*_1$}};
    \filldraw [black] (0.462,0) circle (2pt) node[anchor=north]{\footnotesize{$\sigma_1$}};
    \filldraw [black] (2.163-1,0) circle (2pt) node[anchor=north west]{\footnotesize{$\sigma^*_2$}};
    \filldraw [black] (2.163,0) circle (2pt) node[anchor=north]{\footnotesize{$\sigma_2$}};
    
    \draw[gray,dotted] (0.462,-2) -- (0.462,2);
    \draw[gray,dotted] (2.163,-2) -- (2.163,2);

    \end{tikzpicture}
    \hspace{0.4cm}
     \begin{tikzpicture}[scale=0.75]
    \node (1) at (1.75, 3) {\boxed{\text{\footnotesize{$g(u):=-\arctan(u)$}}}};
    
    \draw[->] (0,0) -- (3.5,0) node[right] {\footnotesize{$\lambda$}};
    \draw[->] (0,-2) -- (0,2);

    \draw [thick, domain=-2:2,variable=\t,smooth,rotate=-90] plot ({(\t)},{(e-1)/(e+1)+(rad(atan((e+1)*\t)))/((e+1)*\t)});
    \draw [thick, domain=-2:2,variable=\t,smooth,rotate=-90] plot ({(\t)},{(e+1)/(e-1)+(rad(atan((e-1)*\t)))/((e-1)*\t)});

    \filldraw [black] (0.462+1,0) circle (2pt) node[anchor=north]{\footnotesize{$\sigma^*_1$}};
    \filldraw [black] (0.462,0) circle (2pt) node[anchor=north]{\footnotesize{$\sigma_1$}};
    \filldraw [black] (2.163+1,0) circle (2pt) node[anchor=north west]{\footnotesize{$\sigma^*_2$}};
    \filldraw [black] (2.163,0) circle (2pt) node[anchor=north]{\footnotesize{$\sigma_2$}};
    
    \draw[gray,dotted] (0.462,-2) -- (0.462,2);
    \draw[gray,dotted] (2.163,-2) -- (2.163,2);

    \end{tikzpicture}  
    \hspace{0.4cm}
     \begin{tikzpicture}[scale=0.75]
    \node (1) at (1.75, 3) {\boxed{\text{\footnotesize{$g(u)=\sqrt{|u|}\sin (u)$}}}};
    
    \draw[->] (0,0) -- (3.25,0) node[right] {\footnotesize{$\lambda$}};
    \draw[->] (0,-2) -- (0,2);

    \draw [thick, domain=-5:5,variable=\t,smooth,rotate=-90] plot ({(0.4*\t)},{(e-1)/(e+1)+(sqrt((e+1)*abs(\t))*(sin(deg((e+1)*\t))))/((e+1)*\t)});

    \draw [thick, domain=-5:5,variable=\t,smooth,rotate=-90] plot ({(0.4*\t)},{(e+1)/(e-1)+(sqrt((e-1)*abs(\t))*(sin(deg((e-1)*\t))))/((e-1)*\t)});




    \filldraw [black] (0.462,0) circle (2pt) node[anchor=north]{\footnotesize{$\sigma_1$}};
    \filldraw [black] (2.163,0) circle (2pt) node[anchor=north]{\footnotesize{$\sigma_2$}};
    
    \draw[gray,dotted] (0.462,-2) -- (0.462,2);
    \draw[gray,dotted] (2.163,-2) -- (2.163,2);

    \end{tikzpicture}  
    \hspace{0.4cm}
     \begin{tikzpicture}[scale=0.75]
    \node (1) at (1.5, 3) {\boxed{\text{\footnotesize{$g(u)=u^2\sin (1/u)$}}}};
    
    \draw[->] (0,0) -- (2.5,0) node[right] {\footnotesize{$\lambda$}};
    \draw[->] (0,-2) -- (0,2);

    \draw [thick, domain=0.001:0.05,variable=\t,smooth,rotate=-90] plot ({(40*\t)},{(e-1)/(e+1)-((e+1)*\t)^(2)*(sin(deg(1/((e+1)*\t))))/((e+1)*\t)});
    \draw [thick, domain=-0.05:-0.001,variable=\t,smooth,rotate=-90] plot ({(40*\t)},{(e-1)/(e+1)-((e+1)*\t)^(2)*(sin(deg(1/((e+1)*\t))))/((e+1)*\t)});

    \draw [thick, domain=0.003:0.05,variable=\t,smooth,rotate=-90] plot ({(40*\t)},{(e+1)/(e-1)-((e-1)*\t)^(2)*(sin(deg(1/((e-1)*\t))))/((e-1)*\t)});
    \draw [thick, domain=-0.05:-0.003,variable=\t,smooth,rotate=-90] plot ({(40*\t)},{(e+1)/(e-1)-((e-1)*\t)^(2)*(sin(deg(1/((e-1)*\t))))/((e-1)*\t)});

    \draw[thick,-] (2.163,-0.15) -- (2.163,0.15);





    \filldraw [black] (0.462,0) circle (2pt) node[anchor=north]{\footnotesize{$\sigma_1$}};
    \filldraw [black] (2.163,0) circle (2pt) node[anchor=north]{\footnotesize{$\sigma_2$}};
    
    \draw[gray,dotted] (0.462,-2) -- (0.462,2);
    \draw[gray,dotted] (2.163,-2) -- (2.163,2);

    \end{tikzpicture}  
    \caption{The bifurcation at infinity for a few examples. \textbf{Left:} There are two (subcritical) bifurcation curves from infinity, each possessing a positive and a negative branch, $\pm c_1(\lambda)$ and $\pm c_2(\lambda)$, which respectively meet the origin at $\sigma^*_1:=\sigma_1-1=-2/(1+\mathrm{e})\approx -0.537\ldots$ and $\sigma^*_2:=\sigma_2-1=2/(\mathrm{e}-1)\approx 1.163\ldots$. \textbf{Middle Left:} The (supercritical) bifurcation curves from infinity meet $0$ at $\sigma^*_1:=\sigma_1+1,\sigma^*_2:=\sigma_2+1$. \textbf{Middle Right:} There are infinitely many turning points accumulating at infinity, see for example~\cite{ARRB10,CP17}. \textbf{Right:} Infinitely many turning points nearby the origin, see \cite{CP12}.}\label{Fig:Bif}
\end{figure}
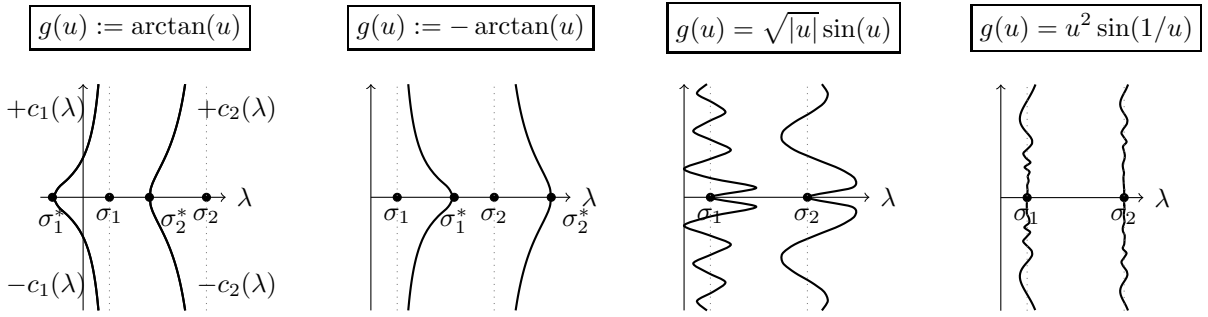

Below we plot the bifurcation curves for some $g(u)$, see Figure~\ref{Fig:Bif}. The main example is $g(u)=\arctan(u)$, which is Lipschitz, odd, monotone, bounded and $g(0)=0,g'(0)=1$. For other $g(u)$ satisfying such properties, the bifurcation diagrams will be qualitatively similar. However, if these properties are violated, 
one would obtain different diagrams.  
Next, we will discuss the hyperbolicity and stability issues of each solution, $u_*(x)\in \{0,\pm u^1(x),\pm u^2(x)\}$, for fixed $\lambda$, which is described according to the eigenvalue problem:
\begin{subequations}\label{EVprob}
\begin{align}
    \varphi_{xx}-\varphi&=\mu(\lambda) \varphi, \qquad\qquad\qquad\qquad x\in (0,1), \label{EVprobEQ}\\    
    -\varphi_x(0)&= \left[\lambda +  g'(u_*(0)) \right]\varphi(0),\label{EVprobbdry1}\\
    \varphi_x(1)&= \left[\lambda +  g'(u_*(1)) \right]\varphi(1).\label{EVprobbdry2}
\end{align}
\end{subequations}
For each fixed $\lambda\in\mathbb{R}$, we will find $\mu=\mu(\lambda)$; we omit the $\lambda$-dependency for simplicity.
First, we discuss the stability of $u_*\equiv 0$, in which $g'(0)=1$. 

For $\mu=-1$, solutions of \eqref{EVprobEQ} are given by $\varphi(x)=ax+\tilde{a}$, whereas the boundary conditions \eqref{EVprobbdry1}-\eqref{EVprobbdry2} imply that $-a=(\lambda+1)\tilde{a}$ and $a=(\lambda+1)(a+\tilde{a})$. Either $\lambda=-1$, yields $a=0$ and consequently $\varphi(x)=\tilde{a}$; or $\lambda\neq -1$, and plugging the first boundary condition into the second implies that $\lambda=1$, and consequently $\varphi(x)=(-2x+1)\tilde{a}$. Therefore, $\mu(\lambda)=-1$ is only an eigenvalue for $\lambda=\pm 1$ with a corresponding eigenfunction. 

For $\mu\neq -1$, general solutions of \eqref{EVprobEQ} are given by $\varphi(x)=a\mathrm{e}^{\sqrt{1+\mu}\, x}+\tilde{a}\mathrm{e}^{-\sqrt{1+\mu}\, x}$ and the boundary conditions \eqref{EVprobbdry1}-\eqref{EVprobbdry2} imply that
\begin{subequations}
\begin{align}
    \tilde{a}&=a \frac{\sqrt{1+\mu}+\lambda+1}{\sqrt{1+\mu}-\lambda-1}\\
    a\mathrm{e}^{\sqrt{1+\mu}}(\sqrt{1+\mu}-\lambda-1)&= \tilde{a}\mathrm{e}^{-\sqrt{1+\mu}}(\sqrt{1+\mu}+\lambda+1).
\end{align}
\end{subequations}
The first equation can be plugged into the second equation, yielding
\begin{equation}\label{identity1}
        \mathrm{e}^{2\sqrt{1+\mu}} = \left(\frac{\sqrt{1+\mu}+\lambda +  1}{\sqrt{1+\mu}-\lambda -  1}\right)^2,
\end{equation}
where the right-hand side has an asymptote at $|1+\mu|=(\lambda+1)^2$. Thus, to solve~\eqref{identity1} for $\mu=\mu(\lambda)$, we compare the graphs of both the left- and right-hand side, where each intersection of both graphs yields an eigenvalue $\mu=\mu(\lambda)$. We emphasize that both sides are equal to 1, when $\mu=-1$.
There are two cases: either $1+\mu>0$ and thereby the square root in~\eqref{identity1} is well-defined, or $1+\mu<0$ and the square root is ill-defined.

\begin{figure}[H]\centering
    \begin{tikzpicture}[scale=0.125]
    \node (1) at (-15, 18) {\boxed{\text{\footnotesize{$\lambda\in \left(\lambda^-_{k},\lambda^-_{k-1}\right), k\neq 0$}}}};
    
    \draw[->] (-30,0) -- (0.1,0) node[right] {\footnotesize{$\mu$}};
    \draw[->] (0,-11) -- (0,11);


    \draw [thick, domain=-(7*pi/2)^2-1+2:-(5*pi/2)^2-1-1.5,variable=\t,smooth] plot ({(0.25*\t)},{tan(deg(sqrt(abs(1+\t))))});
    \draw [thick, domain=-(5*pi/2)^2-1+1.5:-(3*pi/2)^2-1-0.9,variable=\t,smooth] plot ({(0.25*\t)},{tan(deg(sqrt(abs(1+\t))))});
    \draw [thick, domain=-(3*pi/2)^2-1+0.9:-(pi/2)^2-1-0.3,variable=\t,smooth] plot ({(0.25*\t)},{tan(deg(sqrt(abs(1+\t))))});
    \draw [thick, domain=-(pi/2)^2-1+0.3:-1-0.01,variable=\t,smooth] plot ({(0.25*\t)},{tan(deg(sqrt(abs(1+\t))))});    
    
    \draw [lightgray,thick, domain=-37+6.2:0,variable=\t,smooth] plot ({(0.25*\t)},{(2*(-7+1)*(sqrt(abs(1+\t))))/((-7+1)^2-(sqrt(abs(1+\t)))^2)});
    \draw [lightgray,thick, domain=-115:-37-7.5,variable=\t,smooth] plot ({(0.25*\t)},{(2*(-7+1)*(sqrt(abs(1+\t))))/((-7+1)^2-(sqrt(abs(1+\t)))^2)});
    
    \draw[lightgray,dashed] (-0.25*37,-11) -- (-0.25*37,11);

    \filldraw [black] (-1.7,-1) circle (10pt) node[anchor=east]{\footnotesize{$\mu_1$}};
    \filldraw [black] (-6.4,-5.1) circle (10pt) node[anchor=east]{\footnotesize{$\mu_k$}};

    \filldraw [black] (-3.25,-5.1) circle (2pt);
    \filldraw [black] (-4,-5.1) circle (2pt);
    \filldraw [black] (-4.75,-5.1) circle (2pt);

    \filldraw [black] (-14.8,4.2) circle (10pt) node[anchor=east]{\footnotesize{$\mu_{k+1}$}};
    \filldraw [black] (-27.75,1.8) circle (10pt) node[anchor=south east]{\footnotesize{$\mu_{k+2}$}};

    \filldraw [black] (-28.25,1) circle (2pt);
    \filldraw [black] (-29,1) circle (2pt);
    \filldraw [black] (-29.75,1) circle (2pt);

    \end{tikzpicture}
    \hspace{0.5cm}
    \begin{tikzpicture}[scale=0.125]
    \node (1) at (-15, 18) {\boxed{\text{\footnotesize{$\lambda\in \left(\lambda_0^-,-1\right)$}}}};
    
    \draw[->] (-30,0) -- (0.1,0) node[right] {\footnotesize{$\mu$}};
    \draw[->] (0,-11) -- (0,11);


    \draw [thick, domain=-(7*pi/2)^2-1+2:-(5*pi/2)^2-1-1.5,variable=\t,smooth] plot ({(0.25*\t)},{tan(deg(sqrt(abs(1+\t))))});
    \draw [thick, domain=-(5*pi/2)^2-1+1.5:-(3*pi/2)^2-1-0.9,variable=\t,smooth] plot ({(0.25*\t)},{tan(deg(sqrt(abs(1+\t))))});
    \draw [thick, domain=-(3*pi/2)^2-1+0.9:-(pi/2)^2-1-0.3,variable=\t,smooth] plot ({(0.25*\t)},{tan(deg(sqrt(abs(1+\t))))});
    \draw [thick, domain=-(pi/2)^2-1+0.3:-1-0.01,variable=\t,smooth] plot ({(0.25*\t)},{tan(deg(sqrt(abs(1+\t))))});    
    
    \draw [lightgray,thick, domain=-115:-2.4,variable=\t,smooth] plot ({(0.25*\t)},{(2*(-2+1)*(sqrt(abs(1+\t))))/((-2+1)^2-(sqrt(abs(1+\t)))^2)});
    
    \draw[lightgray,dashed] (-0.25*2,-11) -- (-0.25*2,11);

    \filldraw [black] (-0.8,5) circle (10pt) node[anchor=east]{\footnotesize{$\mu_1$}};
    \filldraw [black] (-3.6,0.7) circle (10pt) node[anchor=north]{\footnotesize{$\mu_2$}};
    \filldraw [black] (-11.2,0.5) circle (10pt) node[anchor=north]{\footnotesize{$\mu_{3}$}};
    \filldraw [black] (-24,0.3) circle (10pt) node[anchor=north]{\footnotesize{$\mu_{4}$}};

    \filldraw [black] (-28.25,1) circle (2pt);
    \filldraw [black] (-29,1) circle (2pt);
    \filldraw [black] (-29.75,1) circle (2pt);

    \end{tikzpicture}
    \hspace{0.5cm}
    \begin{tikzpicture}[scale=0.125]
    \node (1) at (-15, 18) {\boxed{\text{\footnotesize{$\lambda=-1$}}}};
    
    \draw[->] (-30,0) -- (0.1,0) node[right] {\footnotesize{$\mu$}};
    \draw[->] (0,-11) -- (0,11);


    \draw [thick, domain=-(7*pi/2)^2-1+2:-(5*pi/2)^2-1-1.5,variable=\t,smooth] plot ({(0.25*\t)},{tan(deg(sqrt(abs(1+\t))))});
    \draw [thick, domain=-(5*pi/2)^2-1+1.5:-(3*pi/2)^2-1-0.9,variable=\t,smooth] plot ({(0.25*\t)},{tan(deg(sqrt(abs(1+\t))))});
    \draw [thick, domain=-(3*pi/2)^2-1+0.9:-(pi/2)^2-1-0.3,variable=\t,smooth] plot ({(0.25*\t)},{tan(deg(sqrt(abs(1+\t))))});
    \draw [thick, domain=-(pi/2)^2-1+0.3:-1-0.01,variable=\t,smooth] plot ({(0.25*\t)},{tan(deg(sqrt(abs(1+\t))))});    
    
    \draw [lightgray,thick, domain=-115:0,variable=\t,smooth] plot ({(0.25*\t)},{(2*(-1+1)*(sqrt(abs(1+\t))))/((-1+1)^2-(sqrt(abs(1+\t)))^2)});
    

    \filldraw [black] (-0.25*1,0) circle (10pt) node[anchor=south west]{\footnotesize{$\mu_1$}};
    \filldraw [black] (-2*pi/2,0) circle (10pt) node[anchor=north]{\footnotesize{$\mu_2$}};
    \filldraw [black] (-6*pi/2,0) circle (10pt) node[anchor=north]{\footnotesize{$\mu_3$}};
    \filldraw [black] (-14*pi/2,0) circle (10pt) node[anchor=north]{\footnotesize{$\mu_{4}$}};

    \filldraw [black] (-28.25,1) circle (2pt);
    \filldraw [black] (-29,1) circle (2pt);
    \filldraw [black] (-29.75,1) circle (2pt);

    \end{tikzpicture}
    \hspace{0.5cm}
    \begin{tikzpicture}[scale=0.125]
    \node (1) at (-15, 18) {\boxed{\text{\footnotesize{$\lambda\in (-1,\lambda^+_0)$}}}};
    
    \draw[->] (-30,0) -- (0.1,0) node[right] {\footnotesize{$\mu$}};
    \draw[->] (0,-11) -- (0,11);


    \draw [thick, domain=-(7*pi/2)^2-1+2:-(5*pi/2)^2-1-1.5,variable=\t,smooth] plot ({(0.25*\t)},{tan(deg(sqrt(abs(1+\t))))});
    \draw [thick, domain=-(5*pi/2)^2-1+1.5:-(3*pi/2)^2-1-0.9,variable=\t,smooth] plot ({(0.25*\t)},{tan(deg(sqrt(abs(1+\t))))});
    \draw [thick, domain=-(3*pi/2)^2-1+0.9:-(pi/2)^2-1-0.3,variable=\t,smooth] plot ({(0.25*\t)},{tan(deg(sqrt(abs(1+\t))))});
    \draw [thick, domain=-(pi/2)^2-1+0.3:-1-0.01,variable=\t,smooth] plot ({(0.25*\t)},{tan(deg(sqrt(abs(1+\t))))});    
    
    \draw [lightgray,thick, domain=-115:-1.16-0.1,variable=\t,smooth] plot ({(0.25*\t)},{(2*(-0.6+1)*(sqrt(abs(1+\t))))/((-0.6+1)^2-(sqrt(abs(1+\t)))^2)});
    
    \draw[lightgray,dashed] (-0.25*1.16,-11) -- (-0.25*1.16,11);

    \filldraw [black] (-0.25*9,-0.4) circle (10pt) node[anchor=north]{\footnotesize{$\mu_2$}};
    \filldraw [black] (-9.5,-0.3) circle (10pt) node[anchor=north]{\footnotesize{$\mu_3$}};


    \filldraw [black] (-22,-0.1) circle (10pt) node[anchor=north]{\footnotesize{$\mu_{4}$}};

    \filldraw [black] (-28.25,-1) circle (2pt);
    \filldraw [black] (-29,-1) circle (2pt);
    \filldraw [black] (-29.75,-1) circle (2pt);

    \end{tikzpicture}
    \hspace{0.5cm}
    \begin{tikzpicture}[scale=0.125]
    \node (1) at (-15, 18) {\boxed{\text{\footnotesize{$\lambda\in \left(\lambda^+_{0},1\right)$}}}};
    
    \draw[->] (-30,0) -- (0.1,0) node[right] {\footnotesize{$\mu$}};
    \draw[->] (0,-11) -- (0,11);


    \draw [thick, domain=-(7*pi/2)^2-1+2:-(5*pi/2)^2-1-1.5,variable=\t,smooth] plot ({(0.25*\t)},{tan(deg(sqrt(abs(1+\t))))});
    \draw [thick, domain=-(5*pi/2)^2-1+1.5:-(3*pi/2)^2-1-0.9,variable=\t,smooth] plot ({(0.25*\t)},{tan(deg(sqrt(abs(1+\t))))});
    \draw [thick, domain=-(3*pi/2)^2-1+0.9:-(pi/2)^2-1-0.3,variable=\t,smooth] plot ({(0.25*\t)},{tan(deg(sqrt(abs(1+\t))))});
    \draw [thick, domain=-(pi/2)^2-1+0.3:-1-0.01,variable=\t,smooth] plot ({(0.25*\t)},{tan(deg(sqrt(abs(1+\t))))});    
    
    \draw [lightgray,thick, domain=-3.24+0.001:-1,variable=\t,smooth] plot ({(0.25*\t)},{(2*(0.8+1)*(sqrt(abs(1+\t))))/((0.8+1)^2-abs(1+\t))});
    \draw [lightgray,thick, domain=-115:-3.24-1.66,variable=\t,smooth] plot ({(0.25*\t)},{(2*(0.8+1)*(sqrt(abs(1+\t))))/((0.8+1)^2-abs(1+\t))});
    
    \draw[lightgray,dashed] (-0.25*3.24,-11) -- (-0.25*3.24,11);
    
    \filldraw [black] (-0.25*3,5) circle (10pt) node[anchor=south east]{\footnotesize{$\mu_2$}};
    \filldraw [black] (-8,-0.8) circle (10pt) node[anchor=north east]{\footnotesize{$\mu_3$}};


    \filldraw [black] (-20.5,-0.4) circle (10pt) node[anchor=north east]{\footnotesize{$\mu_{4}$}};

    \filldraw [black] (-28.25,-1) circle (2pt);
    \filldraw [black] (-29,-1) circle (2pt);
    \filldraw [black] (-29.75,-1) circle (2pt);

    \end{tikzpicture}
    \hspace{0.5cm}
    \begin{tikzpicture}[scale=0.125]
    \node (1) at (-15, 18) {\boxed{\text{\footnotesize{$\lambda=1$}}}};
    
    \draw[->] (-30,0) -- (0.1,0) node[right] {\footnotesize{$\mu$}};
    \draw[->] (0,-11) -- (0,11);


    \draw [thick, domain=-(7*pi/2)^2-1+2:-(5*pi/2)^2-1-1.5,variable=\t,smooth] plot ({(0.25*\t)},{tan(deg(sqrt(abs(1+\t))))});
    \draw [thick, domain=-(5*pi/2)^2-1+1.5:-(3*pi/2)^2-1-0.9,variable=\t,smooth] plot ({(0.25*\t)},{tan(deg(sqrt(abs(1+\t))))});
    \draw [thick, domain=-(3*pi/2)^2-1+0.9:-(pi/2)^2-1-0.3,variable=\t,smooth] plot ({(0.25*\t)},{tan(deg(sqrt(abs(1+\t))))});
    \draw [thick, domain=-(pi/2)^2-1+0.3:-1-0.01,variable=\t,smooth] plot ({(0.25*\t)},{tan(deg(sqrt(abs(1+\t))))});    
    
    \draw [lightgray,thick, domain=-5+0.75:0,variable=\t,smooth] plot ({(0.25*\t)},{(2*(1+1)*(sqrt(abs(1+\t))))/((1+1)^2-(sqrt(abs(1+\t)))^2)});
    \draw [lightgray,thick, domain=-115:-5-0.8,variable=\t,smooth] plot ({(0.25*\t)},{(2*(1+1)*(sqrt(abs(1+\t))))/((1+1)^2-(sqrt(abs(1+\t)))^2)});
    
    \draw[lightgray,dashed] (-0.25*5,-11) -- (-0.25*5,11);

    \filldraw [black] (-0.25*1,0) circle (10pt) node[anchor=south east]{\footnotesize{$\mu_2$}};
    \filldraw [black] (-8,-0.8) circle (10pt) node[anchor=north east]{\footnotesize{$\mu_3$}};


    \filldraw [black] (-20.5,-0.4) circle (10pt) node[anchor=north east]{\footnotesize{$\mu_{4}$}};

    \filldraw [black] (-28.25,-1) circle (2pt);
    \filldraw [black] (-29,-1) circle (2pt);
    \filldraw [black] (-29.75,-1) circle (2pt);

    \end{tikzpicture}
    \hspace{0.5cm}
    \begin{tikzpicture}[scale=0.125]
    \node (1) at (-15, 18) {\boxed{\text{\footnotesize{$\lambda\in \left(1,\lambda^+_{1}\right)$}}}};
    
    \draw[->] (-30,0) -- (0.1,0) node[right] {\footnotesize{$\mu$}};
    \draw[->] (0,-11) -- (0,11);


    \draw [thick, domain=-(7*pi/2)^2-1+2:-(5*pi/2)^2-1-1.5,variable=\t,smooth] plot ({(0.25*\t)},{tan(deg(sqrt(abs(1+\t))))});
    \draw [thick, domain=-(5*pi/2)^2-1+1.5:-(3*pi/2)^2-1-0.9,variable=\t,smooth] plot ({(0.25*\t)},{tan(deg(sqrt(abs(1+\t))))});
    \draw [thick, domain=-(3*pi/2)^2-1+0.9:-(pi/2)^2-1-0.3,variable=\t,smooth] plot ({(0.25*\t)},{tan(deg(sqrt(abs(1+\t))))});
    \draw [thick, domain=-(pi/2)^2-1+0.3:-1-0.01,variable=\t,smooth] plot ({(0.25*\t)},{tan(deg(sqrt(abs(1+\t))))});    
    
    \draw [lightgray,thick, domain=-17+2.8:0,variable=\t,smooth] plot ({(0.25*\t)},{(2*(3+1)*(sqrt(abs(1+\t))))/((3+1)^2-(sqrt(abs(1+\t)))^2)});
    \draw [lightgray,thick, domain=-115:-17-3.3,variable=\t,smooth] plot ({(0.25*\t)},{(2*(3+1)*(sqrt(abs(1+\t))))/((3+1)^2-(sqrt(abs(1+\t)))^2)});
    
    \draw[lightgray,dashed] (-0.25*17,-11) -- (-0.25*17,11);

    \filldraw [black] (-6.4,-4.6) circle (10pt) node[anchor=east]{\footnotesize{$\mu_3$}};    
    \filldraw [black] (-18.7,-1.2) circle (10pt) node[anchor=north east]{\footnotesize{$\mu_{4}$}};

    \filldraw [black] (-28.25,-1) circle (2pt);
    \filldraw [black] (-29,-1) circle (2pt);
    \filldraw [black] (-29.75,-1) circle (2pt);

    \end{tikzpicture}
    \hspace{0.5cm}
    \begin{tikzpicture}[scale=0.125]
    \node (1) at (-15, 18) {\boxed{\text{\footnotesize{$\lambda\in \left(\lambda^+_{k-1},\lambda^+_{k}\right),k\neq 0$}}}};
    
    \draw[->] (-30,0) -- (0.1,0) node[right] {\footnotesize{$\mu$}};
    \draw[->] (0,-11) -- (0,11);


    \draw [thick, domain=-(7*pi/2)^2-1+2:-(5*pi/2)^2-1-1.5,variable=\t,smooth] plot ({(0.25*\t)},{tan(deg(sqrt(abs(1+\t))))});
    \draw [thick, domain=-(5*pi/2)^2-1+1.5:-(3*pi/2)^2-1-0.9,variable=\t,smooth] plot ({(0.25*\t)},{tan(deg(sqrt(abs(1+\t))))});
    \draw [thick, domain=-(3*pi/2)^2-1+0.9:-(pi/2)^2-1-0.3,variable=\t,smooth] plot ({(0.25*\t)},{tan(deg(sqrt(abs(1+\t))))});
    \draw [thick, domain=-(pi/2)^2-1+0.3:-1-0.01,variable=\t,smooth] plot ({(0.25*\t)},{tan(deg(sqrt(abs(1+\t))))});    
    
    \draw [lightgray,thick, domain=-50+9.1:0,variable=\t,smooth] plot ({(0.25*\t)},{(2*(6+1)*(sqrt(abs(1+\t))))/((6+1)^2-(sqrt(abs(1+\t)))^2)});
    \draw [lightgray,thick, domain=-115:-50-10.2,variable=\t,smooth] plot ({(0.25*\t)},{(2*(6+1)*(sqrt(abs(1+\t))))/((6+1)^2-(sqrt(abs(1+\t)))^2)});
    
    \draw[lightgray,dashed] (-0.25*50,-11) -- (-0.25*50,11);

    \filldraw [black] (-4.7,2) circle (10pt) node[anchor=east]{\footnotesize{$\mu_3$}};

    \filldraw [black] (-6.25,5.1) circle (2pt);
    \filldraw [black] (-7,5.1) circle (2pt);
    \filldraw [black] (-7.75,5.1) circle (2pt);
    
    \filldraw [black] (-16.4,-7.5) circle (10pt) node[anchor=east]{\footnotesize{$\mu_{k+2}$}};

    \filldraw [black] (-28.25,-1) circle (2pt);
    \filldraw [black] (-29,-1) circle (2pt);
    \filldraw [black] (-29.75,-1) circle (2pt);

    \end{tikzpicture}
    \caption{Plot of the two graphs in \eqref{identity2} for $\mu<-1$: the (black) left-hand side is given by $\tan(\sqrt{|1+\mu|})$ with asymptotes whenever $\mu=-((2k+1)\pi/2)^2-1,k\in\mathbb{N}_0$, and the (gray) right-hand side is given by $2(\lambda+1)\sqrt{|1+\mu|}/((\lambda+1)^2-|1+\mu|)$ with a (dashed) asymptote at $\mu=-(\lambda+1)^2-1$. The intersection (black dots) of both graphs yields eigenvalues $\mu=\mu(\lambda)<-1$. Note that there are resonances, which occurs whenever both the black and gray curves have asymptote at the same value, which occurs at $\lambda^\pm_k:=\pm (2k+1)\pi/2-1,k\in\mathbb{N}_0$. We emphasize that only two eigenvalues, $\mu_1$ and $\mu_2$, cross the value $\mu=-1$ whenever $\lambda=\pm 1$. 
    }\label{Fig:Bif3}
\end{figure}
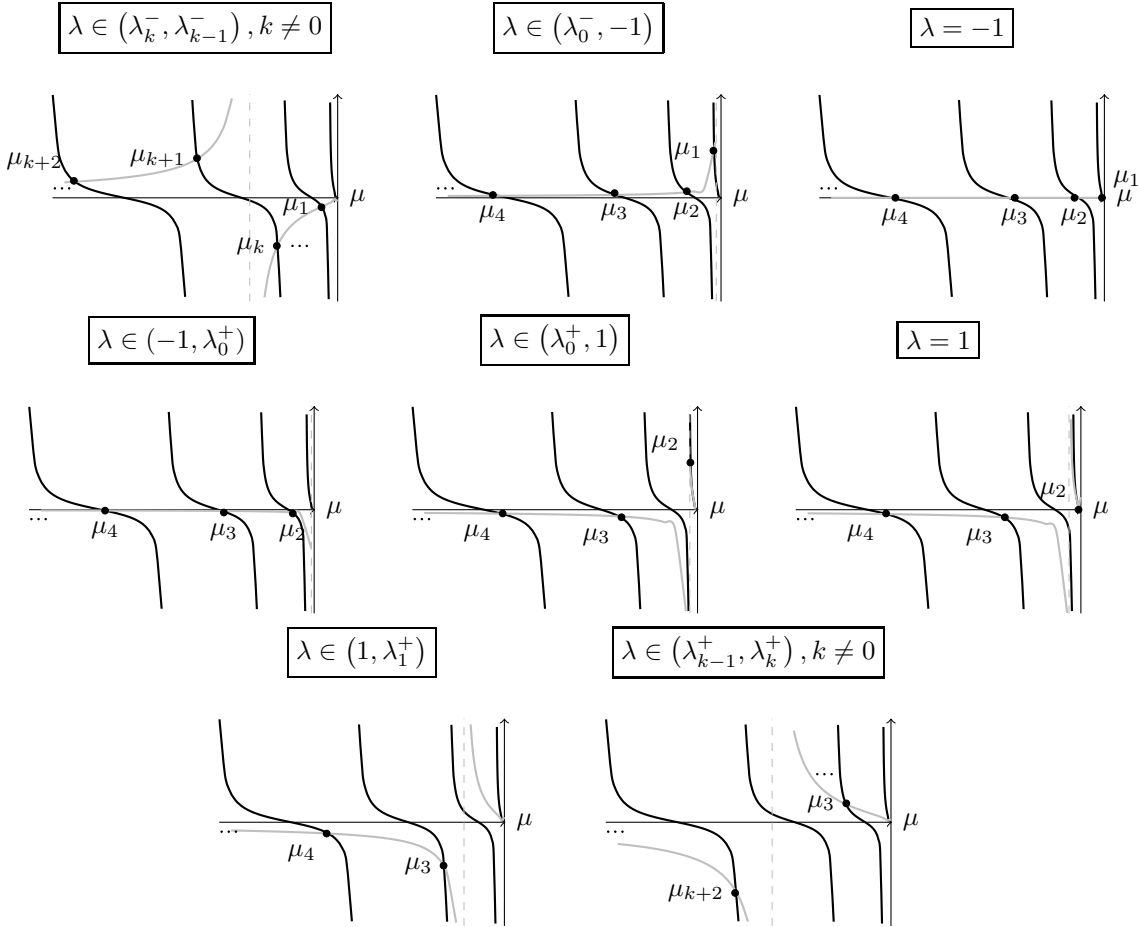

On one hand, if $1+\mu<0$, then the square-root is ill-defined and thus we want to solve
\begin{equation}\label{identity2}
        \mathrm{e}^{2i\sqrt{|1+\mu|}} = \left(\frac{i\sqrt{|1+\mu|}+\lambda +  1}{i\sqrt{|1+\mu|}-\lambda -  1}\right)^2.
\end{equation}
%
%
Comparing the real and imaginary terms of both sides of equation \eqref{identity2}, we obtain
\begin{equation}\label{identity3}
    \tan(\sqrt{|1+\mu|})=\frac{2(\lambda+1)\sqrt{|1+\mu|}}{(\lambda+1)^2-|1+\mu|}.
\end{equation}
We compare the graphs of the left- and right-hand side of the equation above, where each intersection of the graphs yields an eigenvalue $\mu=\mu(\lambda)<-1$. For all $\lambda\in\mathbb{R}$, there are infinitely many such intersections, see Figure~\ref{Fig:Bif3}, and thus infinitely many strictly negative eigenvalues of the trivial solution $u\equiv 0$, which we denote by $\{\mu_n\}_{n=1}^\infty$, with corresponding (stable) eigendirections, which we denote by $\{\varphi_n(x)\}_{n=1}^\infty$. We emphasize that both sides of equation~\eqref{identity3} equal to zero, when $\mu=-1$.
In particular, at $\lambda=1$, the derivative of the left- and right-hand side of \eqref{identity3} coincide at $\mu=-1$.
On the other hand, if $1+\mu>0$, we consider several different subcases, see Figure~\ref{Fig:Bif2}. 
Note that the character of the graphs in \eqref{identity1} change for $\lambda=\pm 1$, as those are the values that yield an eigenvalue $\mu=-1$.
Indeed, if $\lambda<-1$, then the quotient in the right-hand side of \eqref{identity1} is less or equal than 1, whereas the left-hand side is bigger or equal than 1, hence the two graphs do not intersect. 
Monotonicity of the left-hand side of \eqref{identity1} and anti-monotonicity of one of the branches (when $\mu>(\lambda+1)^2-1$) of the right-hand side imply that they always intersect, for $\lambda>-1$.
Moreover, at $\lambda=1$, the derivative of the left- and right-hand side coincide at $\mu=-1$, which occurs when the other branch (when $\mu<(\lambda+1)^2-1$) of the right-hand side crosses the left-hand side, for $\lambda>1$, and therefore a new eigenvalue appears.
This amounts to (at most) two eigenvalues and two corresponding eigenfunctions.
Lastly, the trivial solution bifurcates whenever an eigenvalue cross $\mu=0$. Plugging $\mu=0$ in \eqref{identity1} yields $\mathrm{e}^2=(-(\lambda+2)/\lambda)^2$, which has two solutions, $\lambda =\sigma^*_1,\sigma^*_2$, agreeing with the branches we have found in Figure~\ref{Fig:Bif}.
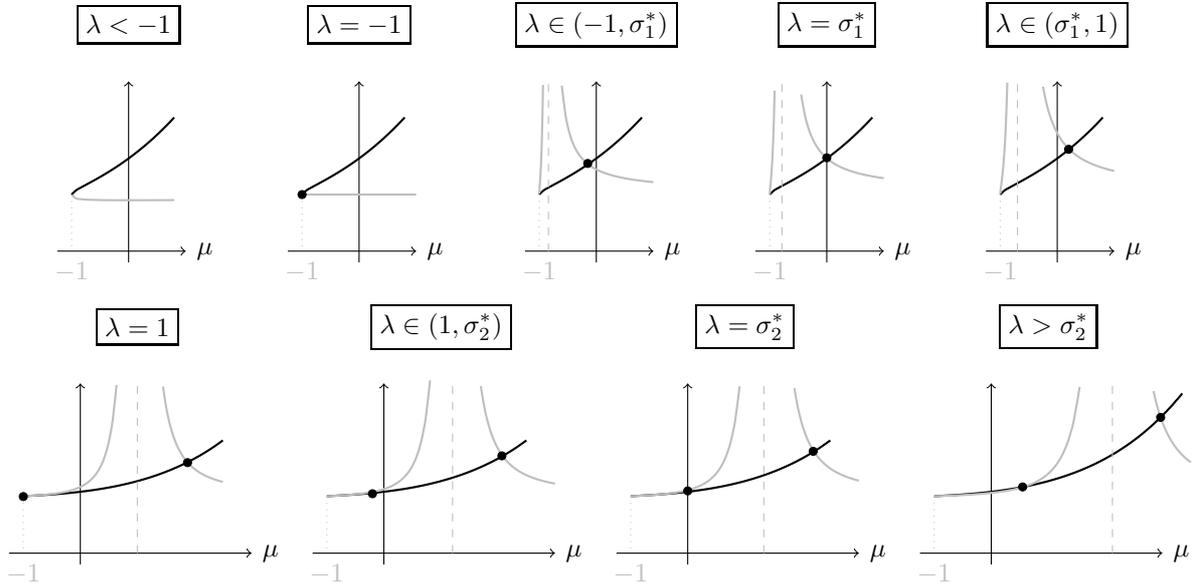
\begin{figure}[H]\centering
    \begin{tikzpicture}[scale=0.75]
    \node (1) at (0, 4) {\boxed{\text{\footnotesize{$\lambda<-1$}}}};
    
    \draw[->] (-1.25,0) -- (1,0) node[right] {\footnotesize{$\mu$}};
    \draw[->] (0,-0.2) -- (0,3);

    
    \draw [lightgray,thick, domain=-1:0.8,variable=\t,smooth] plot ({(\t)},{0.9+((sqrt(1+\t)-2+1)/(sqrt(1+\t)-(-2+1)))^2/10});

    \draw [thick, domain=-1:0.8,variable=\t,smooth] plot ({(\t)},{0.9+exp(2*sqrt(1+\t))/10});

    \draw[lightgray,dotted] (-1,1) -- (-1,0) node[anchor=north]{\footnotesize{$-1$}};

    

    \end{tikzpicture}
    \hspace{0.4cm}
    \begin{tikzpicture}[scale=0.75]
    \node (1) at (0, 4) {\boxed{\text{\footnotesize{$\lambda=-1$}}}};
    
    \draw[->] (-1.25,0) -- (1,0) node[right] {\footnotesize{$\mu$}};
    \draw[->] (0,-0.2) -- (0,3);


    \draw [thick, domain=-1:0.8,variable=\t,smooth] plot ({(\t)},{0.9+exp(2*sqrt(1+\t))/10});
    

    \draw [lightgray,thick, domain=-1:1,variable=\t,smooth] plot ({(\t)},{1});
    
    \draw[lightgray,dotted] (-1,1) -- (-1,0) node[anchor=north]{\footnotesize{$-1$}};


    \filldraw [black] (-1,1) circle (2pt);

    \end{tikzpicture}
    \hspace{0.4cm}
    \begin{tikzpicture}[scale=0.75]
    \node (1) at (0, 4) {\boxed{\text{\footnotesize{$\lambda\in (-1,\sigma^*_1)$}}}};
    
    \draw[->] (-1.25,0) -- (1,0) node[right] {\footnotesize{$\mu$}};
    \draw[->] (0,-0.2) -- (0,3);


    \draw [thick, domain=-1:0.8,variable=\t,smooth] plot ({(\t)},{0.9+exp(2*sqrt(1+\t))/10});
    
    \draw [lightgray,thick, domain=-0.605:1,variable=\t,smooth] plot ({(\t)},{0.9+((sqrt(1+\t)-0.6+1)/(sqrt(1+\t)-(-0.6+1)))^2/10});

    \draw [lightgray,thick, domain=-0.935:-1,variable=\t,smooth] plot ({(\t)},{0.9+((sqrt(1+\t)-0.6+1)/(sqrt(1+\t)-(-0.6+1)))^2/10});
    
    \draw[lightgray,dotted] (-1,1) -- (-1,0) node[anchor=north]{\footnotesize{$-1$}};

    \draw[lightgray,dashed] (0.16-1,0) -- (0.16-1,3);

    \filldraw [black] (-0.15,1.55) circle (2pt);

    \end{tikzpicture}
    \hspace{0.4cm}
    \begin{tikzpicture}[scale=0.75]
    \node (1) at (0, 4) {\boxed{\text{\footnotesize{$\lambda=\sigma^*_1$}}}};
    
    \draw[->] (-1.25,0) -- (1,0) node[right] {\footnotesize{$\mu$}};
    \draw[->] (0,-0.2) -- (0,3);


    \draw [thick, domain=-1:0.8,variable=\t,smooth] plot ({(\t)},{0.9+exp(2*sqrt(1+\t))/10});
    
    \draw [lightgray,thick, domain=-0.45:1,variable=\t,smooth] plot ({(\t)},{0.9+((sqrt(1+\t)-0.537+1)/(sqrt(1+\t)-(-0.537+1)))^2/10});

    \draw [lightgray,thick, domain=-0.915:-1,variable=\t,smooth] plot ({(\t)},{0.9+((sqrt(1+\t)-0.537+1)/(sqrt(1+\t)-(-0.537+1)))^2/10});
    
    \draw[lightgray,dotted] (-1,1) -- (-1,0) node[anchor=north]{\footnotesize{$-1$}};

    \draw[lightgray,dashed] (0.214-1,0) -- (0.214-1,3);

    \filldraw [black] (0,1.65) circle (2pt);

    \end{tikzpicture}
    \hspace{0.4cm}
    \begin{tikzpicture}[scale=0.75]
    \node (1) at (0, 4) {\boxed{\text{\footnotesize{$\lambda\in (\sigma^*_1,1)$}}}};
    
    \draw[->] (-1.25,0) -- (1,0) node[right] {\footnotesize{$\mu$}};
    \draw[->] (0,-0.2) -- (0,3);


    \draw [thick, domain=-1:0.8,variable=\t,smooth] plot ({(\t)},{0.9+exp(2*sqrt(1+\t))/10});
    
    \draw [lightgray,thick, domain=-0.25:1,variable=\t,smooth] plot ({(\t)},{0.9+((sqrt(1+\t)-0.45+1)/(sqrt(1+\t)-(-0.45+1)))^2/10});

    \draw [lightgray,thick, domain=-0.876:-1,variable=\t,smooth] plot ({(\t)},{0.9+((sqrt(1+\t)-0.45+1)/(sqrt(1+\t)-(-0.45+1)))^2/10});
    
    \draw[lightgray,dotted] (-1,1) -- (-1,0) node[anchor=north]{\footnotesize{$-1$}};

    \draw[lightgray,dashed] (0.302-1,0) -- (0.302-1,3);

    \filldraw [black] (0.2,1.8) circle (2pt);

    \end{tikzpicture}
    \hspace{0.4cm}
    \begin{tikzpicture}[scale=0.75]
    \node (1) at (1, 4) {\boxed{\text{\footnotesize{$\lambda=1$}}}};
    
    \draw[->] (-1.25,0) -- (3,0) node[right] {\footnotesize{$\mu$}};
    \draw[->] (0,-0.2) -- (0,3);


    \draw [thick, domain=-1:6,variable=\t,smooth] plot ({(-0.5+0.5*\t)},{1+exp(2*sqrt(1+\t))/200});
    
    \draw [lightgray,thick, domain=3.9:6,variable=\t,smooth] plot ({(-0.5+0.5*\t)},{1+((sqrt(1+\t)+1+1)/(sqrt(1+\t)-(1+1)))^2/200});

    \draw [lightgray,thick, domain=-1:2.27,variable=\t,smooth] plot ({(-0.5+0.5*\t)},{1+((sqrt(1+\t)+1+1)/(sqrt(1+\t)-(1+1)))^2/200});
    
    \draw[lightgray,dotted] (-1,1) -- (-1,0) node[anchor=north]{\footnotesize{$-1$}};

    \draw[lightgray,dashed] (-0.5+0.5*3,0) -- (-0.5+0.5*3,3);

    \filldraw [black] (-1,1) circle (2pt);
    \filldraw [black] (1.88,1.6) circle (2pt);

    \end{tikzpicture}
    %
    %
    \begin{tikzpicture}[scale=0.75]
    \node (1) at (1, 4) {\boxed{\text{\footnotesize{$\lambda\in (1,\sigma^*_2)$}}}};
    
    \draw[->] (-1.25,0) -- (3,0) node[right] {\footnotesize{$\mu$}};
    \draw[->] (0,-0.2) -- (0,3);


    \draw [thick, domain=-1:6,variable=\t,smooth] plot ({(-0.5+0.5*\t)},{1+exp(2*sqrt(1+\t))/200});
    
    \draw [lightgray,thick, domain=4.4:7,variable=\t,smooth] plot ({(-0.5+0.5*\t)},{1+((sqrt(1+\t)+1.1+1)/(sqrt(1+\t)-(1.1+1)))^2/200});

    \draw [lightgray,thick, domain=-1:2.62,variable=\t,smooth] plot ({(-0.5+0.5*\t)},{1+((sqrt(1+\t)+1.1+1)/(sqrt(1+\t)-(1.1+1)))^2/200});
    
    \draw[lightgray,dotted] (-1,1) -- (-1,0) node[anchor=north]{\footnotesize{$-1$}};

    \draw[lightgray,dashed] (-0.5+0.5*3.41,0) -- (-0.5+0.5*3.41,3);

    \filldraw [black] (2.07,1.72) circle (2pt);

    \filldraw [black] (-0.2,1.05) circle (2pt);
    \end{tikzpicture}
    %
    %
    \begin{tikzpicture}[scale=0.75]
    \node (1) at (1, 4) {\boxed{\text{\footnotesize{$\lambda =\sigma^*_2$}}}};
    
    \draw[->] (-1.25,0) -- (3,0) node[right] {\footnotesize{$\mu$}};
    \draw[->] (0,-0.2) -- (0,3);


    \draw [thick, domain=-1:6,variable=\t,smooth] plot ({(-0.5+0.5*\t)},{1+exp(2*sqrt(1+\t))/200});
    
    \draw [lightgray,thick, domain=4.74:7,variable=\t,smooth] plot ({(-0.5+0.5*\t)},{1+((sqrt(1+\t)+1.163+1)/(sqrt(1+\t)-(1.163+1)))^2/200});

    \draw [lightgray,thick, domain=-1:2.82,variable=\t,smooth] plot ({(-0.5+0.5*\t)},{1+((sqrt(1+\t)+1.163+1)/(sqrt(1+\t)-(1.163+1)))^2/200});
    
    \draw[lightgray,dotted] (-1,1) -- (-1,0) node[anchor=north]{\footnotesize{$-1$}};

    \draw[lightgray,dashed] (-0.5+0.5*3.67,0) -- (-0.5+0.5*3.67,3);

    \filldraw [black] (2.2,1.8) circle (2pt);

    \filldraw [black] (0,1.1) circle (2pt);

    \end{tikzpicture}
    %
    %
    \begin{tikzpicture}[scale=0.75]
    \node (1) at (1, 4) {\boxed{\text{\footnotesize{$\lambda >\sigma^*_2$}}}};
    
    \draw[->] (-1.25,0) -- (3,0) node[right] {\footnotesize{$\mu$}};
    \draw[->] (0,-0.2) -- (0,3);


    \draw [thick, domain=-1:7.7,variable=\t,smooth] plot ({(-0.5+0.5*\t)},{1+exp(2*sqrt(1+\t))/200});
    
    \draw [lightgray,thick, domain=6.69:8,variable=\t,smooth] plot ({(-0.5+0.5*\t)},{1+((sqrt(1+\t)+1.5+1)/(sqrt(1+\t)-(1.5+1)))^2/200});

    \draw [lightgray,thick, domain=-1:4.1,variable=\t,smooth] plot ({(-0.5+0.5*\t)},{1+((sqrt(1+\t)+1.5+1)/(sqrt(1+\t)-(1.5+1)))^2/200});
    
    \draw[lightgray,dotted] (-1,1) -- (-1,0) node[anchor=north]{\footnotesize{$-1$}};

    \draw[lightgray,dashed] (-0.5+0.5*5.25,0) -- (-0.5+0.5*5.25,3);

    \filldraw [black] (2.97,2.4) circle (2pt);

    \filldraw [black] (0.55,1.17) circle (2pt);

    \end{tikzpicture}
    \caption{Plot of the two graphs in \eqref{identity1} for $\mu>-1$: the (black) left-hand side, $\mathrm{e}^{2\sqrt{1+\mu}}$, and the (gray) right-hand side, $(\sqrt{1+\mu}+\lambda+1)^2/(\sqrt{1+\mu}-\lambda-1)^2$, with a (dashed) asymptote at $\mu=(\lambda+1)^2-1$. Each intersection (black dots) of both graphs amounts to an eigenvalue $\mu=\mu(\lambda)>-1$. We emphasize that there are at most two eigenvalues that cross $\mu=0$, which occur exactly at the values $\lambda=\sigma^*_1,\sigma^*_2$, and thereby yield the bifurcation of the trivial equilibrium. }\label{Fig:Bif2}
\end{figure}
Therefore we have found all eigenvalues $\mu=\mu(\lambda)$ and their corresponding eigenfunctions of the linearized operator at the trivial equilibrium, $u_*\equiv 0$. 
We have shown that there are infinitely many eigenvalues $\{\mu_n\}_{n=1}^\infty$ satisfying
\begin{equation}\label{EV}
    \ldots < \mu_{n+1} < \mu_n < \mu _{n-1} < \ldots < \mu_3 < \mu_2 < \mu_1,
\end{equation}
with respective eigenfunctions $\{\varphi_n(x)\}_{n=1}^\infty$. 
On one hand, $\mu_n<0$ for all $n\leq 3$, which correspond to stable eigendirections. On the other hand, 
the two eigenvalues $\mu_1(\lambda),\mu_2(\lambda)$ cross the value $\mu=0$ whenever $\lambda=\sigma^*_1,\sigma^*_2$, and correspond to unstable eigendirections, see Figure~\ref{Fig:Bif2}. 

\begin{figure}[H]\centering
    \begin{tikzpicture}[scale=0.75]
    \node (1) at (0, 4) {\boxed{\text{\footnotesize{$\lambda\in (\sigma^*_1,\sigma_1)$}}}};
    
    \draw[->] (-1.25,0) -- (1,0) node[right] {\footnotesize{$\mu$}};
    \draw[->] (0,-0.2) -- (0,3);

    \draw [thick, domain=-1:0.8,variable=\t,smooth] plot ({(\t)},{0.9+exp(2*sqrt(1+\t))/10});

    \draw [lightgray,thick, domain=-1:0.8,variable=\t,smooth] plot ({(\t)},{0.9+((sqrt(1+\t)+(-0.3+0.460444))/(sqrt(1+\t)-(-0.3+0.460444)))^2/10});



    \draw[lightgray,dotted] (-1,1) -- (-1,0) node[anchor=north]{\footnotesize{$-1$}};

    \filldraw [black] (-0.69,1.21) circle (2pt) node[anchor=north]{\footnotesize{$\mu_1$}};
    
    \end{tikzpicture}
    \hspace{3cm}
    \begin{tikzpicture}[scale=0.7]
    \node (1) at (1, 4) {\boxed{\text{\footnotesize{$\lambda\in (\sigma^*_2,\sigma_2)$}}}};
    
    \draw[->] (-1.25,0) -- (3,0) node[right] {\footnotesize{$\mu$}};
    \draw[->] (0,-0.2) -- (0,3);

    \draw [thick, domain=-1:5.5,variable=\t,smooth] plot ({(-0.5+0.5*\t)},{1+exp(2*sqrt(1+\t))/100});

    \draw [lightgray,thick, domain=-1:2.25,variable=\t,smooth] plot ({(-0.5+0.5*\t)},{1+((sqrt(1+\t)+(1.2+0.89642))/(sqrt(1+\t)-(1.2+0.89642)))^2/100});

    \draw [lightgray,thick, domain=4.95:6.5,variable=\t,smooth] plot ({(-0.5+0.5*\t)},{1+((sqrt(1+\t)+(1.2+0.89642))/(sqrt(1+\t)-(1.2+0.89642)))^2/100});

    \draw[lightgray,dotted] (-1,1) -- (-1,0) node[anchor=north]{\footnotesize{$-1$}};

    \draw[lightgray,dashed] (-0.5+0.5*3.84,2.9) -- (-0.5+0.5*3.84,0) node[anchor=north]{\footnotesize{$\mu=(\lambda+1)^2-1$}};

    \filldraw [black] (-0.65,1.05) circle (2pt) node[anchor=north]{\footnotesize{$\mu_2$}};

    \filldraw [black] (2.075,2.42) circle (2pt) node[anchor=west]{\footnotesize{$\mu_1$}};

    \end{tikzpicture}
    \caption{The graphs in \eqref{identity1} for $\mu>-1$: $\mathrm{e}^{2\sqrt{1+\mu}}$ in black and $(\sqrt{1+\mu}+\lambda+g'(u_*(0)))^2/(\sqrt{1+\mu}-\lambda-g'(u_*(0)))^2$ in gray, where $u_*\in \{\pm u^1,\pm u^2\}$. 
    \textbf{Left:} For $u_*=\pm u^1$ with $\lambda \in (\sigma^*_1,\sigma_1)$, there is one intersections of both graphs yielding an eigenvalue $\mu_1=\mu_1(\lambda)\in (-1,0)$, and hence $\pm u^1$ are stable hyperbolic equilibra, since all other eigenvalues satisfy $\mu_k<-1,k=2,3,\ldots$. Moreover, $\mu_1 \to 0$ if either $\lambda\to\sigma^*_1$ or $\sigma_1$. 
    \textbf{Right:} For $u_*=\pm u^2$ with $\lambda \in (\sigma^*_2,\sigma_2)$, there are two intersections of both graphs yielding the eigenvalues $\mu_1=\mu_1(\lambda)\in (0,\infty)$ and $\mu_2=\mu_2(\lambda)\in [-1,0)$, and thus $\pm u^1$ are hyperbolic saddles, since all other eigenvalues satisfy $\mu_k<-1,k=3,4,\ldots$. Moreover, $\mu_2 \to 0$ if either $\lambda\to\sigma^*_2$ or $\sigma_2$. 
    }\label{Fig:Bif4}
\end{figure}
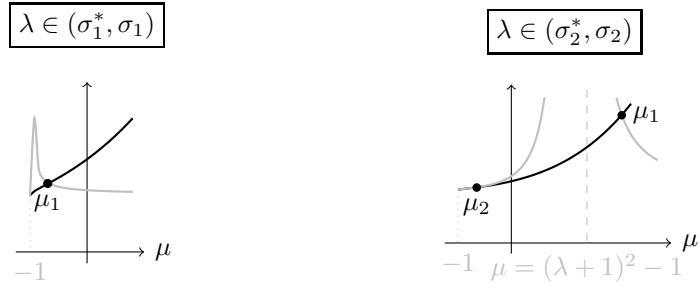

Similarly, we can prove the hyperbolicity and stability of the non-trivial equilibria, $u_*(x)\in \{\pm u^1(x),\pm u^2(x)\}$. 
Indeed, the equilibria $\pm u^1(x)$ are stable for $\lambda \in (\sigma^*_1,\sigma_1)$, whereas $\pm u^2(x)$ are saddles for $\lambda \in (\sigma^*_2,\sigma_2)$. See Figure~\ref{Fig:Bif4}.
%
%
%
Hence, 
the above discussion can be summarized as follows.
\begin{lem} \label{lem:eq} 
Under the above conditions on $g(u)$, all equilibria $u_*\equiv0,\pm u^1(x), \pm u^2(x)$ of \eqref{intro:PDE} are hyperbolic for $\lambda\neq \sigma^*_1,\sigma_1,\sigma^*_2,\sigma_2$. Moreover, the following holds:
\begin{enumerate} 
    \item[(i)] If $\lambda\in (-\infty,\sigma^*_1)$, there is a unique equilibrium, $u_*\equiv 0$, which is asymptotically stable. At $\lambda=\sigma^*_1$, the solution $u_*\equiv 0$ undergoes a supercritical pitchfork bifurcation.
    \item[(ii)] If $\lambda\in (\sigma^*_1,\sigma_1)$, then there are three equilibria: $u_*\equiv 0$ is one-dimensional unstable and $\pm u^1 (x)$ are asymptotically stable. As $\lambda\to \sigma_1$, then $\pm u^1(x)$ escape to infinity towards $\pm\Phi_1(x)$ as in~\eqref{EQconvSteklov}.
    \item[(iii)] If $\lambda\in (\sigma_1,\sigma^*_2)$, then $u_*\equiv 0$ is the only bounded equilibrium, which is one-dimensional unstable.  At $\lambda=\sigma^*_2$, there is another supercritical pitchfork bifurcation of $u_*\equiv 0$. 
    \item[(iv)] If $\lambda\in (\sigma^*_2,\sigma_2)$, then there are three equilibria: $u_*\equiv 0$ is two-dimensional unstable and $\pm u^2 (x)$ are saddles with one-dimensional unstable direction. As $\lambda\to \sigma_2$, then $\pm u^2(x)$ escape to infinity towards $\pm\Phi_2(x)$ as in~\eqref{EQconvSteklov}.
    \item[(v)] If $\lambda\in (\sigma_2,\infty)$, the only equilibrium is $u_*\equiv 0$, which is two-dimensional unstable.
\end{enumerate}
\end{lem}
There is another linear eigenvalue problem which plays a role in the global dynamics. It corresponds to a linearization of \eqref{intro:PDE} ``at infinity'', in which  $\lim_{||u||\to\infty} g'(u)\to 0$: 
\begin{subequations}\label{EVprobinfty}
\begin{align}
    \varphi_{xx}-\varphi&=\mu^\infty(\lambda) \varphi, \qquad\qquad\qquad x\in (0,1), \label{EVprobEQinfty}\\    
    -\varphi_x(0)&= \lambda\varphi(0),\label{EVprobbdry1infty}\\
    \varphi_x(1)&= \lambda \varphi(1).\label{EVprobbdry2infty}
\end{align}
\end{subequations}
Note that the eigenvalues are $\mu^\infty_n(\lambda)=\mu_n(\lambda-1),n\in\mathbb{N}$, where $\mu_n(\lambda-1)$ are the eigenvalues of the linearization at $u_*\equiv 0$ in \eqref{EVprob}, whenever $g'(0)=1$. Therefore, the eigenvalues ``at infinity'' correspond to a translation of the eigenvalues at $u_*\equiv 0$ by a factor $g'(0)=1$, however, the eigenfunctions $\varphi_n(x)$ are the same.

Next, we discuss the growth of general time-dependent solutions in each projected eigendirection $\varphi_n(x)$. This will also enable us to discuss dynamical issues, such as the existence of blow-up solutions and their corresponding rescaled asymptotics. In particular, we will obtain the existence of bounded and unbounded global attractors.

\subsection{Asymptotic behavior}\label{sec:att}
 
To describe the global dynamics of solutions of the parabolic problem \eqref{intro:PDE}, we will project the evolution of solutions of \eqref{intro:PDE}, given by the variation of constants formula in~\eqref{VoC}, into each eigendirection of the orthonormal basis $\{\varphi_n(x)\}_{n=1}^\infty$ of $L^2$. 
This amounts to:
\begin{equation}
    u(t,x)=\sum_{n=1}^\infty u_n(t) \varphi_n(x) ,
\end{equation}
where we define $u_n(t):=\langle u(t,x),\varphi_n(x) \rangle_{L^2}$.

We will show that $u_n(t),n\geq 3,$ are bounded as $t\to \infty$, whereas $u_1(t),u_2(t)$ may blow up in infinite time, respectively, when $\lambda>\sigma_1,\sigma_2$, with an accurate asymptotic description. 
%

%
%
\begin{lem} \label{lem:nondiss} Under the above conditions on $g(u)$, 
the following holds:
\begin{enumerate} 
    \item[(i)] If $\lambda<\sigma_{1}$, then 
    there is a compact global attractor, $\mathcal{A}_\lambda$. If $\lambda>\sigma_{1}$, then 
    there is an unbounded attractor, $\mathcal{A}_\lambda$, containing some infinite-time blow-up solutions. 
    \item[(ii)] If $\lambda\in (\sigma_{1},\sigma_2)$, then any blow-up solution satisfies
            \begin{equation}\label{BupSigma1}
            \lim_{t\to \infty} \frac{u(t,x)}{||u(t,x)||_{L^2}}\to \iota \varphi_1 , \text{ for some } \iota =-1,+1.
        \end{equation}
        In particular, there are initial data in $W^u(u_*\equiv 0)$ that blow-up.
        Moreover, for $\lambda\in (\sigma_2^*,\sigma_2)$, there are also initial data in $W^u(\pm u^2(x))$ that blow-up.
    \item[(iii)] if $\lambda>\sigma_2$, then any blow-up solution satisfies
        \begin{equation}\label{BupSigma2}
            \lim_{t\to \infty} \frac{u(t,x)}{||u(t,x)||_{L^2}}\to \iota \varphi_N , \text{ for some } \iota= -1,+1 \text{ and some } N=1,2.
        \end{equation}
        In particular, there are initial data in $W^u(u_*\equiv 0)$ that blow-up as in \eqref{BupSigma2} for each $\iota=+1,-1$ and $N=1,2$.
\end{enumerate}
\end{lem}
\begin{pf}
The result for $\lambda<\sigma_{1}$ has been achieved in~\cite[Section 7]{ARRB07}. 
Next, for $\lambda>\sigma_{1}$, we prove which eigenmodes $u_n(t)$ remain bounded and which will possibly grow -- recall that solutions are global in time.
From~\eqref{intro:PDE2}, each $u_n(t)$ satisfy the projected equation 
\begin{equation}
   \frac{d}{dt}u_{n}(t) = \langle u_{t}(t,x),\varphi_{n}(x)\rangle_{L^{2}} = 
\langle u_{xx}-u,\varphi_{n}(x)\rangle_{L^{2}}+\langle G(u(t,x)),\varphi_{n}(x)\rangle_{L^{2}} \ , 
\end{equation}
which, by selfadjointness of $A:=\partial_x^2 - Id$ and $G_{n}(t):=\langle G(u(t,x)),\varphi_{n}(x)\rangle_{L^{2}}$, implies
\begin{equation}\label{eq:ODE}
    \frac{d}{dt}u_{n}(t) = \mu_n^\infty(\lambda)u_{n}(t)+ G_{n}(t). 
\end{equation} 
%
%
A general solution of \eqref{eq:ODE} is given by the variation of constants formula:
\begin{eqnarray}\label{eq:FormulaDaSolucao}
    u_{n}(t)=u_{n}(0)e^{\mu_n^\infty(\lambda)t}+\int_{0}^{t}e^{\mu_n^\infty(\lambda)(t-s)}G_{n}(s)ds \ ,
\end{eqnarray} 
where $u_{n}(0)=\langle u_{0}(x),\varphi_{n}(x)\rangle_{L^{2}}$. Thus, if $\mu_n^\infty(\lambda)<0$, we obtain 
    $|u_{n}(t)| \leq |u_{n}(0)| + \frac{C}{\mu_n^\infty(\lambda)} $,
where $C>0$ is an uniform bound for $G(u)$. 
Indeed, note that the boundedness of $g(u)$ implies that $G(u)$ is bounded 
and thereby $G_{n}(t)$ is uniformly bounded for all $n\in \mathbb{N}$. 

Note that if $\lambda<\sigma_1$, then $\mu_n^\infty(\lambda)<0$ for all $n\in\mathbb{N}$, and hence all $u_n(t)$ remain bounded in time.
If $\lambda\in (\sigma_{1},\sigma_2)$, then $\mu_n^\infty(\lambda)<0$ for all $n>1$ and $\mu_1^\infty(\lambda)>0$; and hence all $u_n(t)$ remain bounded for $n>1$, but we can not conclude that $u_1(t)$ is bounded.
If $\lambda>\sigma_2$, then $\mu_n^\infty(\lambda)<0$ for all $n>2$ and $\mu_1^\infty(\lambda),\mu_2^\infty(\lambda)>0$; and hence all $u_n(t)$ remain bounded for $n>2$, but we can not conclude that $u_1(t),u_2(t)$ are bounded.

Next, we analyze the eigenmodes $u_n(t)$ corresponding to $\mu_n^\infty(\lambda)>0$ and show that they amount to growing solutions. Note that the integral in the right hand side of \eqref{eq:FormulaDaSolucao} is unbounded as $t\to\infty$, since $\mu_n^\infty(\lambda)(t-s)>0$, and thus we rewrite \eqref{eq:FormulaDaSolucao} as
\begin{equation}\label{eq:particular}
    u_{n}(t)=u_{n}^{h}(0)e^{\mu^\infty_n(\lambda)t}+\int_{\infty}^{t}e^{\mu^\infty_n(\lambda)(t-s)}G_{n}(s)ds \, ,
\end{equation} 
where $u_{n}^{h}(0) := u_{n}(0) - \int_{\infty}^{0}e^{\mu_n(\lambda)s}G_{n}(s)ds$.\footnote{Note that a general solution can be regarded in the form 
$u_{n}(t)=u_{n}^{h}(t)+u_{n}^{p}(t)$, where $u_{n}^{h}(t)$ is a solution of the corresponding homogeneous equation and $u_{n}^{p}(t)$ is a particular 
solution of \eqref{eq:ODE}.  We take the particular solution  $u_{n}^{p}(t)=\int_{\infty}^{t}e^{\mu_n(\lambda)(t-s)}G_{n}(s)ds$. }
In this case, since $G_{n}(t)$ is uniformly bounded for all $n\in \mathbb{N}$ and $\mu^\infty_n(\lambda)(t-s)<0$, the integral of the right hand side of \eqref{eq:particular} remains bounded. Therefore the growth of $u_n(t)$ is dictated by $u_{n}^{h}(0)e^{\mu^\infty_n(\lambda)t}$.
Hence, there is some $u_{n}^{h}(0)\neq 0$, which in turn implies that $u_{n}^{h}(t)=u_{n}^{h}(0)e^{\mu_n(\lambda)t}$ grows exponentially as $t\to \infty$ and thus there is a solution that blows up in infinite time. This finishes the proof of the second part of (i), which claims that there are infinite-time blow-up solutions and thereby the system is non-dissipative. To obtain the unbounded attractor, we apply~\cite{FernandesBortolan}.

To prove $(ii)$, we use the $L^2$-inner product and orthonormality of $\{\varphi_n(x)\}_{n\in\mathbb{N}}$ to obtain:
\begin{align}
    \lim_{t\rightarrow\infty}\left|\left|\frac{u(t,\cdot)}{\| u(t,\cdot)\|_{L^2}}-\iota \varphi_{N}(\cdot)\right|\right|^2_{L^2} 
    & = 2 - 2\lim_{t\rightarrow\infty}\frac{\langle u(t,\cdot),\iota\varphi_{N}(\cdot)\rangle_{L^2}}{\|u(t,\cdot)\|_{L^2}} \\
& = 2 - \lim_{t\rightarrow\infty}2\iota\frac{u_{N}(t)}{( \sum_{n=1}^{\infty}u_{n}^{2}(t))^{\frac12}} = 0
\end{align}
if, and only if,
\begin{equation}
   \lim_{t\rightarrow\infty}\frac{u_{N}(t)}{( \sum_{n=1}^{\infty}u_{n}^{2}(t))^{\frac12}} = \textcolor{black}{\iota} \ . 
\end{equation}
If $\lambda\in (\sigma_1,\sigma_2)$, note that the only growing mode is $u_{1}(t)$, whereas all other $u_{n}(t)$ with $n>1$ are bounded. Therefore, 
\begin{equation}
    \lim_{t\rightarrow\infty}\frac{u_1(t)}{( \sum_{n=1}^{\infty}u_{n}^{2}(t))^{\frac12}} = \textcolor{black}{\iota} \qquad \text{ and } \qquad   \lim_{t\rightarrow\infty}\frac{u_m(t)}{( \sum_{n=1}^{\infty}u_{n}^{2}(t))^{\frac12}} = 0,
\end{equation} 
for all $m>1$. 

For $\lambda\in (\sigma_1,\sigma_2^*)$, we know that $u_*\equiv 0$ has a one-dimensional unstable manifold, due to Lemma~\ref{lem:eq}. An initial data in this unstable manifold satisfying $u_1^h(0)\neq 0$ gives the desired blow-up solution.
Similarly, we choose an appropriate initial data in the unstable manifold of $\pm u^2(x)$, in case that $\lambda\in (\sigma_2^*,\sigma_2)$. That concludes the proof of $(ii)$.

If $\lambda>\sigma_2$, the only possible growing modes are $u_{1}(t),u_{2}(t)$, whereas all other $u_{n}(t)$ with $n>2$ are necessarily bounded. There are two cases. First, if $u_0(x)$ is such that $u_1^h(0)= 0$, then $u_1(t)$ remains bounded and $u_2(t)$ is the unique growing mode, and we have that
\begin{equation}
    \lim_{t\rightarrow\infty}\frac{u_2(t)}{( \sum_{n=1}^{\infty}u_{n}^{2}(t))^{\frac12}} = \textcolor{black}{\iota} \qquad \text{ and } \qquad   \lim_{t\rightarrow\infty}\frac{u_m(t)}{( \sum_{n=1}^{\infty}u_{n}^{2}(t))^{\frac12}} = 0,
\end{equation} 
for all $m\neq 2$. Second, if $u_1^h(0)\neq 0$, then 
\begin{equation}
    \lim_{t\rightarrow\infty}\frac{u_1(t)}{( \sum_{n=1}^{\infty}u_{n}^{2}(t))^{\frac12}} = \textcolor{black}{\iota} \qquad \text{ and } \qquad   \lim_{t\rightarrow\infty}\frac{u_m(t)}{( \sum_{n=1}^{\infty}u_{n}^{2}(t))^{\frac12}} = 0
\end{equation} 
for all $m>1$. 
Note that the growth to infinity of $u_{1}(t)$ is determined by the term $u_{1}^{h}(0)e^{\mu_1(\lambda)t}$, whereas the growth of $u_2(t)$ is determined by $u_{2}^{h}(0)e^{\mu_2(\lambda)t}$, if $u_2^h(0)\neq 0$. Thus, $u_1(t)$ grows faster than $u_{2}(t)$, since $\mu_1(\lambda)>\mu_2(\lambda)$.

For any $u_0(x)$ that yields a solution which blows up in infinite time, we conclude that the rescaled solution $\frac{u(t,x)}{\| u(t,x) \|_{L^2}}$ converges to either $\varphi_N$ or $-\varphi_N$, with $N=1,2$. From a zero number stabilization argument and the sign condition $\iota$ on $u(t,0)$, for all sufficiently large $t \in (T,\infty)$, we get the limit to $\iota \varphi_N$.

For $\lambda>\sigma_2$, we know that $u_*\equiv 0$ has a two-dimensional unstable manifold, due to Lemma~\ref{lem:eq}. As before, an initial data in this unstable manifold satisfying $u_1^h(0)\neq 0$ or $u_2^h(0)\neq 0$, provides the desired blow-up solutions. That concludes the proof of $(iii)$.
\qed
\end{pf}
We analyze the structure of the attractor in~\eqref{char:A}, $\mathcal{A}_\lambda$, under the above hypothesis on $g(u)$. 
If $\lambda<\sigma^*_1$, the attractor consists of the equilibrium $u_*\equiv 0$.
At $\lambda=\sigma^*_1$, such solution undergoes a supercritical pitchfork bifurcation and emanates two equilibria $\pm u^1(x)$ that become arbitrarily large as $\lambda\nearrow\sigma_1$. However, $\mathcal{A}_\lambda$ remains bounded for fixed $\lambda\in (\sigma^*_1,\sigma_1)$ and consists of the equilibria $u_*\equiv 0,\pm u^1(x)$ and their heteroclinics. 
Note that their heteroclinics persist for all such parameters, due to hyperbolicity and transversality of the stable and unstable manifolds; see \cite[Section 4]{RochaFiedler} for a proof that can be adapted to this case.
If $\lambda>\sigma_1$, 
the attractor $\mathcal{A}_\lambda$ is unbounded and contains the blow-up solutions within the unstable manifold of bounded equilibria in Lemma~\ref{lem:nondiss}, due to~\eqref{char:A}. 
%

Next, we provide a `compactified' perspective of such unbounded behavior. We will describe the limiting objects of such heteroclinics to infinity using a Poincaré projection, which will be equilibria of the projected semiflow at infinity. This justifies such nomenclature of \emph{heteroclinics to infinity}. In particular, the Steklov eigenfunctions, $\pm \Phi_1(x),\pm \Phi_2(x)$, are equilibria at infinity for $\lambda=\sigma_1,\sigma_2$, as they are the respective limit of the unbounded bifurcation curves of equilibria $\pm u^1(x),\pm u^2(x)$ as $\lambda\to \sigma_1,\sigma_2$, in Section~\ref{sec:ubdd}.

\subsection{Poincaré Projection }\label{sec:P} 
%
%
The infinite-time blow-up solutions were interpreted as heteroclinics to infinity in \cite{Hell09}. In order to describe the dynamics of unbounded solutions, an infinite dimensional sphere $\mathcal{S}^\infty$ was added at infinity with an induced semiflow, by means of a Poincaré projection. 
Previous investigations to understand such structure at infinity for semilinear equations were done for $f(u)$ in \cite{Hell09,BenGal10}, for $f(x,u,u_x)$ in \cite{RochaPimentel16}, and  periodic boundary condition in \cite{P}. 
We now describe this process in detail.

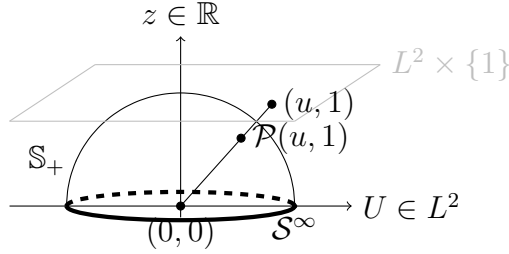
\begin{figure}[H]\centering
    \begin{tikzpicture}[scale=1.5]
    \draw[->] (-1.5,0) -- (1.5,0) node[right] {$U\in L^2 $};
    \draw[->] (0,-0.1) -- (0,1.5) node[above] {$z\in \mathbb{R}$};

    \draw [domain=0:3.14,variable=\t,smooth] plot ({cos(\t r)},{sin(\t r)}); 

    \draw[ultra thick] (-1,0) arc (180:360:1cm and 0.125cm);
    \draw[ultra thick,dashed] (-1,0) arc (180:0:1cm and 0.125cm);

    \draw[-] (0,0) -- (0.8,0.9);    


    \draw[lightgray,-] (-0.75,1.25) -- (1.75,1.25) node[right] {$L^2\times \{ 1\}$};
    \draw[lightgray,-] (-1.5,0.75) -- (1,0.75);
    \draw[lightgray] (-0.75,1.25) -- (-1.5,0.75); 
    \draw[lightgray] (1.75,1.25) -- (1,0.75);

    \filldraw [black] (0.8,0.9) circle (1pt) node[anchor=west]{$(u,1)$};
    \filldraw [black] (0,0) circle (1pt) node[anchor=north]{$(0,0)$};
    \filldraw [black] (0.53,0.6) circle (1pt) node[anchor=west]{$\mathcal{P}(u,1)$};
    \filldraw [black] (-0.915,0.4) circle (0.01pt) node[anchor=east]{$\mathbb{S}_+$};
    \filldraw [black] (1,0) circle (0.5pt) node[anchor=north]{$\mathcal{S}^\infty$};
        
    \end{tikzpicture}
\caption{Poincaré projection $\mathcal{P}$ from phase-space $X\hookrightarrow L^2\times \{ 1\}$ into the hemisphere $\mathbb{S}_+$ through collinearity with $(0,0)$. As solutions of equation \eqref{intro:PDE} blow-up in infinite time, $||u(t)||_{L^2}\to \infty$, the projection $\mathcal{P}(u,1)$ converges to the equator $\mathcal{S}^\infty:=\mathbb{S}_+|_{z=0}$.}\label{Fig:P}
\end{figure}

The Poincaré projection maps the phase-space $X$ of \eqref{intro:PDE} to a subset of the unit sphere in $L^2\times \mathbb{R}$. 
Indeed, identify $X$ with $X\times\{1\}\subseteq L^2\times\{1\}$. The set $L^2\times\{1\}$ is the tangent space (at the north pole) of the northern hemisphere, $\mathbb{S}_+\subseteq L^2\times \mathbb{R}$, called the \emph{Poincaré hemisphere} and given by
\begin{equation}\label{Phemis}
    \mathbb{S}_+:=\{(U,z)\in L^2\times [0,1]: ||U||^{2}_{L^2}+z^{2}=1\}.
\end{equation}
Then for each point in phase-space, $u\in X\subseteq L^2$, consider the line that passes through $(u,1)\in L^2\times [0,1]$ and the origin $(0,0)\in L^2\times [0,1]$. This line intersects the upper hemisphere $\mathbb{S}_+$ at a point, which defines the projection $\mathcal{P}:X\times\{1\}\to \mathbb{S}_+$, called the \emph{Poincaré projection}. See Figure~\ref{Fig:P}.
%
The coordinates of the projection $\mathcal{P}(u,1)$ are denoted by $(U,z)$ and can be computed from the colinearity of the points $(0,0)$, $(u,1)$, and the intersection at $(U,z)$ with the hemisphere $\mathbb{S}_+$, yielding 
\begin{equation}\label{P}
    (U,z):={\mathcal{P}}(u,1)=\left(\frac{u}{ \sqrt{1+||u||_{L^2}^2}},\frac{1}{ \sqrt{1+||u||_{L^2}^2}}\right).
\end{equation}
%
Therefore, Hell's perspective that infinite-time blow-up solutions in $X$ are heteroclinics to infinity can be interpreted (in terms of the Poincaré projection $\mathcal{P}$) as heteroclinics in the hemisphere $\mathbb{S}_+$ that converge to the equator $\mathcal{S}^\infty$, which is characterized by $z=0$ and $||U||^2_{L^2}=1$. Indeed, note that $z=1$ if, and only if, $u \equiv 0$; hence the origin of $X$ is mapped to the north pole of $\mathbb{S}_+$. Similarly, $z$ decreases to $0$ if, and only if, $||u||_{L^2}$ increases to $\infty$; 
Thus the asymptotic dynamics of the projected unbounded semiflow is contained in the unit sphere of $L^2$, consisting of bounded trajectories with coordinates $(U,0)$. For this reason, the equator of $\mathbb{S}_+$ is called the \textit{sphere at infinity}, and it is denoted by $\mathcal{S}^\infty:=\mathbb{S}_+|_{z=0}$.

The projection $\mathcal{P}$ induces a semiflow on $\mathbb{S}_+|_{z>0}$, which is obtained by a homothety (with scale factor $z$) of the original vector field in \eqref{intro:PDE}. 
%
Indeed, differentiating \eqref{P} with respect to time, the new variables $(U,z)$ satisfy
\begin{subequations}\label{flowSPHERE}
\begin{align}
    U_t&
    = U_{xx}-U- \langle U_{xx}-U,U\rangle_{L^2} \cdot U, \label{flowSPHERE1}\\ 
    -U_x(0)&
    = \lambda U (0)+g^z(U(0)) ,\label{bdry1}\\
    U_x(1)&= \lambda U (1)+g^z(U(1)) ,\label{bdry2}\\
    z_t& 
    =-\langle U_{xx}-U,U\rangle_{L^2} \cdot z,\label{flowSPHERE2}
\end{align}
\end{subequations}
where $g^z(U):=z g(z^{-1}U)$ is a homothety of the nonlinearity in the boundary conditions with the scale factor $z>0$. 
The projection \eqref{P} thereby induces a semiflow within $\mathbb{S}_+|_{z>0}$ described by equations \eqref{flowSPHERE}, which consists of an ODE for $z$ coupled to a nonlocal semilinear parabolic PDE for $U$, since the inner product can be rewritten as:
\begin{equation}\label{innprod}
    -\langle U_{xx}-U,U\rangle_{L^2} = -\left (U_xU\right)|_{x=0}^{x=1} + ||U||_{H^1}.
\end{equation}
%
%
%
Moreover, the induced (nonlinear and nonlocal) semiflow at $\mathcal{S}^\infty$, i.e. for $z=0$, is given by the limit as $z\to 0$, which may produce, in a number of settings, a degenerate or singular semiflow at the sphere at infinity $\mathcal{S}^\infty$; see \cite{BenGal10,Hell09} for further details and examples.
Note that in the present case, such limit is well-defined, since the inner product in~\eqref{innprod} is bounded and therefore we can compute the limit in equation~\eqref{flowSPHERE2} as $z\to 0$, which implies that $\mathcal{S}^\infty$ is invariant.
The semiflow at $\mathcal{S}^\infty$, as $z\to 0$, is given by
\begin{subequations}\label{flowSPHEREatinfty}
\begin{align}
    U_t&
    = U_{xx}-U-\langle U_{xx}-U,U\rangle_{L^2} \cdot U,  \label{flowSPHERE1infty}\\ 
    -U_x(0)&
    = \lambda U (0), \label{bdry1infty}\\
    U_x(1)&= \lambda U (1), \label{bdry2infty}
\end{align}
\end{subequations}
since $g$ is bounded and hence $g^z(U)\to 0$ as $z\to 0$. Note that in the limit, the inner product becomes $-\langle U_{xx}-U,U\rangle_{L^2} = - \lambda \left (U^2(1) + U^2(0)\right) + ||U||_{H^1}.$

Therefore, the Steklov eigenfunctions $\Phi_1(x),\Phi_2(x)$, which satisfy $0=U_{xx}-U$ and the boundary conditions~\eqref{bdry1infty}-\eqref{bdry2infty}, are equilibria at infinity when $\lambda=\sigma_1,\sigma_2$. Therefore, through the compactification in~\eqref{P}, not only we obtain a limiting object at infinity (corresponding to a Steklov eigenfunction), recovering the known results of \cite{ARRB07,ARRB09} for unidimensional spatial domain, but more importantly, that such limiting object are equilibria of the induced semiflow at infinity.

More generally, all the eigenfunctions $\{\varphi_n\}_{n\in\mathbb{N}}$ of the linearization at infinity, given by equation~\eqref{EVprobinfty}, are solutions of \eqref{flowSPHEREatinfty} and thereby equilibria at infinity for all $\lambda\in\mathbb{R}$, given in coordinates as
\begin{eqnarray}
    \varphi_n^\infty:=(\varphi_{n},0)\in\mathcal{S}^\infty \quad \text{ for all } n\in\mathbb{N}.
\end{eqnarray}
In the particular case that $\lambda=\sigma_1,\sigma_2$, we have $\varphi_1(x)=\Phi_1(x),\varphi_2(x)=\Phi_2(x)$.
However, due to Lemma~\ref{lem:nondiss}, only the eigenfunctions $\pm\varphi_1(x),\pm\varphi_2(x)$ can be reached by infinite-time blow-up solutions. 
Next, we seek to describe the full dynamics of the induced compactified semiflow at the sphere at infinity. 
\begin{lem}
    There are four solutions of \eqref{flowSPHEREatinfty} which are heteroclinics orbits at the sphere at infinity from $\pm\varphi^\infty_2(x)$ to $\pm\varphi^\infty_1(x)$ as $t$ increases. 
\end{lem}
\begin{pf}
To account for the growth in each eigendirection, we consider two further projections, $\mathcal{P}_{N}$, for each fixed $N=1,2$, from $X\times \{1\}$ into hyperplanes $C_{N}$, which are tangent to the equilibria at infinity, $(\varphi_{N},0)\in\mathcal{S}^\infty$, and given by
\begin{equation}
    C_{N}:=\{ (U,z)\in L^2\times [0,1]\text{ $|$ } U_{N}=+1,U_n\in\mathbb{R}\text{ for all } n\in\mathbb{N} \backslash\{N\}\},
\end{equation}
where  $U_n:=\langle U,\varphi_{n}\rangle$ for all $n\in\mathbb{N}$.
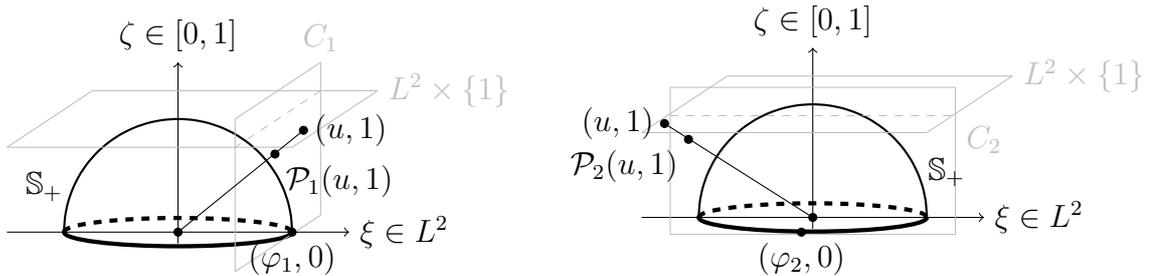
\begin{figure}[H]\centering
    \hspace{0.75cm}
    \begin{tikzpicture}[scale=1.5]
    \draw[->] (-1.5,0) -- (1.5,0) node[right] {$\xi\in L^2 $};
    \draw[->] (0,-0.1) -- (0,1.5) node[above] {$\zeta\in [0,1]$};

    \draw [thick, domain=0:3.14,variable=\t,smooth] plot ({cos(\t r)},{sin(\t r)}); 

    \draw[ultra thick] (-1,0) arc (180:360:1cm and 0.125cm);
    \draw[ultra thick,dashed] (-1,0) arc (180:0:1cm and 0.125cm);
    

    \draw[lightgray,-] (-0.75,1.25) -- (1.75,1.25) node[right] {$L^2\times \{ 1\}$};
    \draw[lightgray,-] (-1.5,0.75) -- (1,0.75);
    \draw[lightgray] (-0.75,1.25) -- (-1.5,0.75); 
    \draw[lightgray] (1.75,1.25) -- (1,0.75); 

    \draw[lightgray] (1.25,1.5) -- (0.5,1); 
    \draw[lightgray] (1.25,0.15) -- (0.5,-0.35); 
    \draw[lightgray] (0.5,-0.35) -- (0.5,1); 
    \draw[lightgray] (1.25,0.15) -- (1.25,1.5)  node[above] {$C_{1}$};     
    \draw[lightgray,dashed] (1.25,1.25) -- (0.5,0.75); 
    
    
    \filldraw [black] (-0.915,0.4) circle (0.01pt) node[anchor=east]{$\mathbb{S}_+$};
    \filldraw [black] (1.1,0.9) circle (1pt) node[anchor=west]{$(u,1)$};
    \filldraw [black] (0,0) circle (1pt); 
    \filldraw [black] (0.85,0.69) circle (1pt) node[anchor=north west]{$\mathcal{P}_{1}(u,1)$};
    \filldraw [black] (1,0) circle (1pt) node[anchor=north]{$(\varphi_{1},0)$};

    \draw[-] (0,0) -- (1.1,0.9);
    
    \end{tikzpicture}
    \hspace{0.25cm}
    \begin{tikzpicture}[scale=1.5]
    \draw[->] (-1.5,0) -- (1.5,0) node[right] {$\xi\in L^2 $};
    \draw[->] (0,-0.1) -- (0,1.5) node[above] {$\zeta\in [0,1]$};

    \draw [thick, domain=0:3.14,variable=\t,smooth] plot ({cos(\t r)},{sin(\t r)}); 

    \draw[ultra thick] (-1,0) arc (180:360:1cm and 0.125cm);
    \draw[ultra thick,dashed] (-1,0) arc (180:0:1cm and 0.125cm);
    

    \draw[lightgray,-] (-0.75,1.25) -- (1.75,1.25) node[right] {$L^2\times \{ 1\}$};
    \draw[lightgray,-] (-1.5,0.75) -- (1,0.75);
    \draw[lightgray] (-0.75,1.25) -- (-1.5,0.75); 
    \draw[lightgray] (1.75,1.25) -- (1,0.75); 

    \draw[lightgray] (-1.25,1.15) -- (1.25,1.15); 
    \draw[lightgray] (-1.25,-0.15) -- (1.25,-0.15); 
    \draw[lightgray] (-1.25,-0.15) -- (-1.25,1.15); 
    \draw[lightgray] (1.25,-0.15) -- (1.25,1.15);     
    \draw[lightgray,dashed] (-1.25,0.9) -- (1.25,0.9) node[anchor=north west] {$C_{2}$};   
    
    \filldraw [black] (0.915,0.4) circle (0.01pt) node[anchor=west]{$\mathbb{S}_+$};
    \filldraw [black] (-1.3,0.83) circle (1pt) node[anchor=east]{$(u,1)$};
    \filldraw [black] (0,0) circle (1pt); 
    \filldraw [black] (-1.09,0.69) circle (1pt) node[anchor=north east]{$\mathcal{P}_{2}(u,1)$};
    \filldraw [black] (-0.1,-0.13) circle (1pt) node[anchor=north]{$(\varphi_{2},0)$};

    \draw[-] (0,0) -- (-1.3,0.83);
    
    \end{tikzpicture}
\caption{Projections $\mathcal{P}_{N}$ for each $N=1,2$ from phase-space $X\hookrightarrow L^2\times \{ 1\}$ into the hyperplanes $C_{N}$, which are respectively tangent (at the point $(\varphi_{N},0)$) to the sphere at infinity $\mathcal{S}^\infty$.}\label{fig:Ptilde}
\end{figure}
%
%
Consider any point $u\in X\subseteq L^2$ and a line that passes through the points $(u,1),(0,0)\in  L^2\times [0,1]$. The intersection of such line with the plane $C_{N}$ defines the projection $\mathcal{P}_{N}$, see Figure \ref{fig:Ptilde}, with coordinates given by $(\xi,\zeta)$:
\begin{equation}\label{Ptilde}
    (\xi,\zeta):=\mathcal{P}_{N}(u,1)=\left(\frac{u}{ \langle u,\varphi_{N} \rangle_{L^2}},\frac{1}{ \langle u,\varphi_{N} \rangle_{L^2}}\right).
\end{equation}
We emphasize that we abuse the notation and we use the same symbol $(\xi,\zeta)$ given by~\eqref{Ptilde} to denote \emph{different} variables for each of the planes $C_N,N=1,2$. 
Note the plane $C_{N}$ can be rewritten in the coordinates $(\xi,\zeta)$ as
\begin{equation}
    C_{N}:=\{ (\xi,\zeta)\in L^2\times \mathbb{R}\text{ $|$ } \xi_{N}=+1,\xi_n\in\mathbb{R}\text{ for all } n\in\mathbb{N} \backslash\{N\}\},
\end{equation}
where  $\xi_n:=\langle\xi,\varphi_{n}\rangle_{L^2}$ for all $n\in\mathbb{N}$.

The projection $\mathcal{P}_{N}$ induces a semiflow on $C_{N}$ through differentiation of \eqref{Ptilde} with respect to time, yielding 
\begin{subequations}\label{flowCj*}
\begin{align}
    \xi_t&
    =\xi_{xx} - \xi-\langle \xi_{xx} - \xi,\varphi_{N} \rangle_{L^2} \cdot \xi  \label{xit}\\
    -\xi_x(0)&
    = \lambda\xi (0)+g^\zeta(\xi(0)) ,\label{bdry1xi}\\
    \xi_x(1)&= \lambda\xi (1)+g^\zeta(\xi(1)) ,\label{bdry2xi}\\
    \zeta_t&
    =-\langle \xi_{xx} - \xi, \varphi_{N} \rangle_{L^2} \cdot \zeta. \label{zetat}
\end{align}
\end{subequations}
where the projected vector field is a homothety of the original vector field \eqref{intro:PDE} with scale factor $\zeta:=\langle u,\varphi_{N} \rangle^{-1}_{L^2}$, i.e. $g^\zeta(\xi):=\zeta g(\zeta^{-1}\xi)$, as in \eqref{flowSPHERE}.
Note that the inner product $\langle \xi_{xx} - \xi, \varphi_{N} \rangle_{L^2} $ is bounded, since $\langle\xi,\varphi_{N}\rangle_{L^2}=1$ in $C_{N}$, and thus the equation~\eqref{zetat} is well-defined and we can compute the limit as $\zeta\to 0$, since $g^\zeta(\xi)\to 0$:
\begin{subequations}\label{flowCj*zeta0}
\begin{align}
    \xi_t&
    =\xi_{xx} - \xi-\langle \xi_{xx} - \xi,\varphi_{N} \rangle_{L^2} \cdot \xi  \label{xitzeta0}\\
    -\xi_x(0)&
    = \lambda\xi (0) ,\label{bdry1xizeta0}\\
    \xi_x(1)&= \lambda\xi (1) ,\label{bdry2xizeta0}\\
    \zeta_t&
    =0. \label{zetatzeta0}
\end{align}
\end{subequations}
Note that the coordinates $(U,z)$ and $(\xi,\zeta)$ are related by means of colinearity,
\begin{equation}\label{changeofcoord}
    (\xi,\zeta)=\frac{1}{\langle U , \varphi_{N}\rangle_{L^2}}(U,z),
\end{equation}
and thereby the semiflow described by both variables are locally topologically conjugate. In other words, the planes $C_N$ provide suitable coordinates to describe the evolution of solutions in the hemisphere $\mathbb{S}_+$.

In order to dissect the induced semiflow of equation \eqref{xit} in $C_{N}|_{\zeta=0}$, we write each $\varphi_n$-eigendirectional semiflow as $\xi_n:=\langle \xi ,\varphi_n\rangle_{L^2}$, in the limit $\zeta\to 0$, which satisfies
\begin{align}\label{xijlim}
   (\xi_n)_t &= \langle \xi_{xx}-\xi,\varphi_n \rangle_{L^2}- \langle  \xi_{xx}-\xi,\varphi_{N}  \rangle_{L^2} \cdot \xi_n\nonumber\\ 
   &=(\mu_{n}-\mu_{{N}})\xi_n,
\end{align}
since $A=\partial_{xx}-Id$ is self-adjoint and $\xi_N=+1$ at the hyperplane $C_N$. Note that since each $\varphi_n$ is an eigenfunction of the eigenvalue problem, it satisfies the linear Robin boundary conditions, and thereby so does $\xi(t,x)=\sum_n \xi_n(t) \varphi_n(x)$.

Therefore, infinite-time blow-up behavior of the solutions $u(t)$ induces the linear flow \eqref{xijlim} in the projected coordinates $(\xi,0)\in C_{N}|_{\zeta=0}$. In particular, equilibria of \eqref{xijlim} occur when $\xi_{n}=0$ for all $n\in\mathbb{N}$, except $\xi_{N}\neq 0$. The only of these equilibria that lie on the sphere at infinity $\mathcal{S}^\infty$ is when $\xi_{N}=\pm 1$. Thus, the equilibria at infinity are given by the eigenfunctions, $\pm\varphi^\infty_{N}:=(\pm\varphi_{N},0)\in\mathcal{S}^\infty$, and can be characterized as
\begin{align}\label{EQcsi}
    \pm\varphi^\infty_{N} &:= \{ (\xi,0)\in \mathcal{S}^\infty:\xi_{N}=\pm 1, \mbox{ and } \xi_{n}=0 \ \forall n \neq N \},\\
    &=\{ (U,0)\in \mathcal{S}^\infty:U_{N}=\pm 1, \mbox{ and } U_{n}=0 \ \forall n \neq N \} \ .
\end{align}
%
Let us now construct the heteroclinic network at infinity, $\mathcal{H}^\infty\subseteq \mathcal{S}^\infty$. 
Note that, due to Lemma~\ref{lem:nondiss}, only the equilibria $\varphi^\infty_N$ with $N=1,2$ are attainable from infinite-time blow-up solutions. To compactify the unbounded attractor, we will show that there is a heteroclinic connection from $\pm\varphi^\infty_{2}$ to $\pm\varphi^\infty_{1}$ whenever $\lambda>\sigma_2$. 

Consider the hyperplane $C_{2}$, where $\zeta=0$, which contains the equilibrium $\pm\varphi^\infty_{2}$.
Then, each coordinate $\xi_n(t)$ belongs to a linear subspace where its evolution is given by \eqref{xijlim}. 
When $\lambda>\sigma_2$, we have that the eigenvalues satisfy $\mu_n<\ldots<\mu_3<0<\mu_2<\mu_1$. 
Thus, we have the following:

\quad \boxed{\text{$n=1$}} \, $\xi_1'=(\mu_1-\mu_2)\xi_1$, where $\mu_1-\mu_2>0$ and thus $\xi_1(t)$ increases.

\quad \boxed{\text{$n=2$}} \, $\xi_2'=0$ and thus $\xi_2\equiv +1$.    

\quad \boxed{\text{$n\geq 3$}} \, $\xi_n'=(\mu_n-\mu_2)\xi_n$, where $\mu_n-\mu_2<0$ and thus $\xi_n(t)$ decreases.    

Considering the change of coordinates in~\eqref{changeofcoord} in the plane $C_2$, note that
\begin{equation}\label{hetinfty}
        (\xi_1(t),1,\xi_3(t),\ldots)=\xi=\frac{1}{\langle U , \varphi_{2}\rangle_{L^2}} U = \left( \frac{U_1(t)}{U_2(t)},  1, \frac{U_3(t)}{U_2(t)},\ldots\right),
\end{equation}
According to the projection $\mathcal{P}_2$, the coordinate $\xi_1(t)\to\infty$ as $t\to \infty$, which implies that $\frac{U_1(t)}{U_2(t)}\to \infty$ as $t\to \infty$, i.e., $U_1(t)$ grows faster than $U_2(t)$. On the other hand, $\xi_n(t)\to 0$ as $t\to \infty$ for $n\geq 3$.

Consider now the hyperplane $C_{1}$, where $\zeta=0$, which contains the equilibrium $\pm\varphi^\infty_{1}$.
Then, each coordinate $\xi_n(t)$ belongs to a linear subspace where its evolution is given by \eqref{xijlim}, and for $\lambda>\sigma_2$ we obtain 

    \quad \boxed{\text{$n=1$}} \, $\xi_1'=0$ and thus $\xi_1\equiv +1$.  
    
    \quad \boxed{\text{$n\geq 2$}} \, $\xi_n'=(\mu_n-\mu_1)\xi_1$, where $\mu_n-\mu_1<0$ and thus $\xi_n(t)$ decreases.    

Considering the change of coordinates in~\eqref{changeofcoord} in the plane $C_1$, note that
\begin{equation}\label{hetinfty2}
        (\xi_1(t),1,\xi_3(t),\ldots)=\xi=\frac{1}{\langle U , \varphi_{1}\rangle_{L^2}} U = \left( 1, \frac{U_2(t)}{U_1(t)},  \frac{U_3(t)}{U_1(t)},\ldots\right),
\end{equation}
Since $U_n(t)$ are the coordinates of the solution in the hemisphere at infinity that leaves from $\varphi_2$, according to~\eqref{hetinfty}, we know that $\frac{U_1(t)}{U_2(t)}\to \infty$ as $t\to \infty$, we have that $\frac{U_2(t)}{U_1(t)}\to 0$ and thereby $\xi_2(t)\to 0$. Similarly for the others $\xi_n(t)\to 0,n\geq 3$.
Thus, we obtain the heteroclinic network in the sphere at infinity, $\mathcal{S}^\infty$. 
\qed
\end{pf}
%

\subsection{Main Theorem}\label{sec:Thm}

\begin{thm} \label{mainthm}
Under the above hypothesis on $g(u)$, there are five distinct dynamical regimes for the global attractor of equation \eqref{intro:PDE} as $\lambda$ crosses each of the bifurcating eigenvalues $\sigma^*_1,\sigma_2,\sigma^*_2,\sigma_2$, as follows; see Figure~\ref{Fig:P1d}.
\begin{enumerate} 
    \item[(i)] For $\lambda\in (0,\sigma^*_1)$, all solutions remain bounded and $u_*\equiv 0$ is the global attractor. 
    \item[(ii)] For $\lambda\in (\sigma^*_1,\sigma_1)$, solutions are bounded and the global attractor consists of the three equilibria, $0$ and $\pm u^1 (x)$, and two heteroclinics from $0$ to $\pm u^1 (x)$ as $t$ increases. Moreover, as $\lambda\to \sigma_1$, then $\pm u^1 (x)$ escape to infinity towards $\pm\Phi_1(x)$, as in \eqref{EQconvSteklov}. 
    \item[(iii)] For $\lambda\in (\sigma_1,\sigma^*_2)$, the compactified unbounded global attractor consists of $u_*\equiv 0$, two equilibria at infinity, $\pm\varphi^\infty_1(x)$, and two unbounded heteroclinics (i.e. infinite-time blow-up solutions) from $0$ to $\pm \varphi^\infty_1 (x)$ as $t$ increases, as in \eqref{BupSigma1}. 
    \item[(iv)] For $\lambda\in (\sigma^*_2,\sigma_2)$, the compactified unbounded attractor consists of three bounded equilibria, $ 0$ and $\pm u^2(x)$, two equilibria at infinity, $\pm\varphi^\infty_1(x)$, and their heteroclinics. 
    There are bounded heteroclinics from $0$ to $\pm u^2 (x)$ and unbounded heteroclinics from both $0$ and $\pm u^2(x)$ to $\pm\varphi^\infty_1(x)$, as $t$ increases, as in \eqref{BupSigma1}.
    Moreover, as $\lambda\to \sigma_2$, then $\pm u^2 (x)$ escape to infinity towards $\pm\Phi_2(x)$, as in \eqref{EQconvSteklov}. 
    \item[(v)] For $\lambda\in (\sigma_2,\infty)$, the compactified unbounded global attractor consists of $u_*\equiv 0$, four equilibria at infinity, $\pm\varphi^\infty_1(x),\pm\varphi^\infty_2(x)$, and their heteroclinics. 
    There are unbounded heteroclinics from $u_*\equiv 0$ to both $\pm\varphi^\infty_1(x)$ and $\pm\varphi^\infty_2(x)$, as $t$ increases, as in \eqref{BupSigma2}. At infinity, there are heteroclinics from $\pm\varphi^\infty_2(x)$ to $\pm\varphi^\infty_1(x)$, as $t$ increases. 
\end{enumerate}
\end{thm}
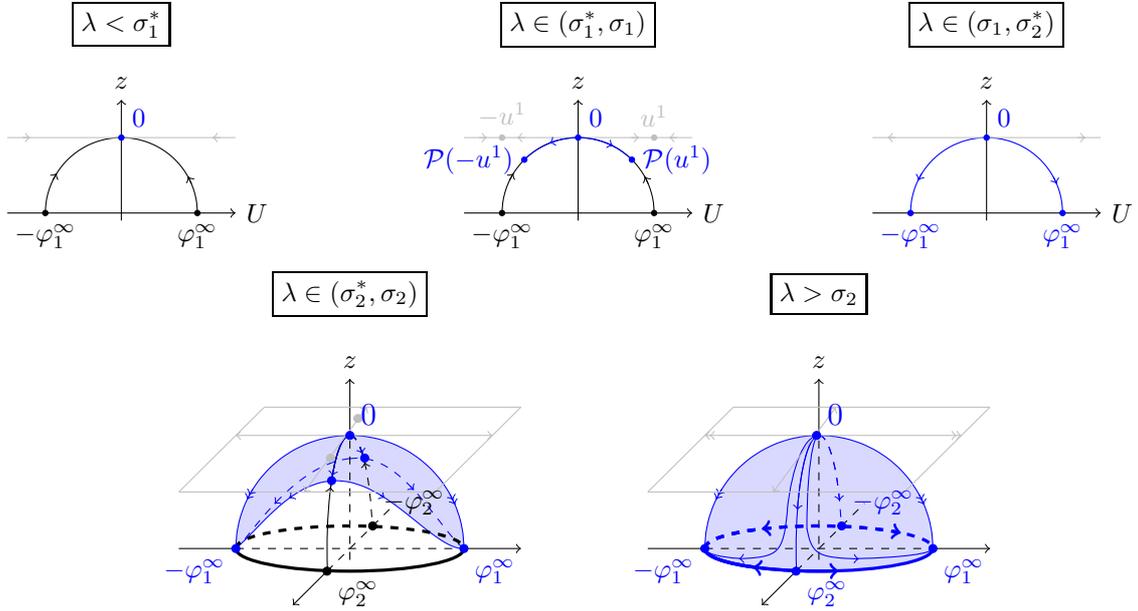
\begin{figure}[H]\centering
    \begin{tikzpicture}[scale=1]
    \node (1) at (0, 2.5) {\boxed{\text{\footnotesize{$\lambda<\sigma^*_1$}}}};    

    \draw[->] (-1.5,0) -- (1.5,0) node[right] {\footnotesize{$U$}};
    \draw[->] (0,-0.1) -- (0,1.5) node[above] {\footnotesize{$z$}};

    \draw [domain=0:3.14,variable=\t,smooth] plot ({cos(\t r)},{sin(\t r)}); 
    
    \draw[lightgray,-] (-1.5,1) -- (1.5,1);

    \filldraw [blue] (0,1) circle (1pt) node[anchor=south west]{\footnotesize{$0$}};
    \filldraw [black] (1,0) circle (1pt) node[anchor=north]{\footnotesize{$\varphi_1^\infty$}};
    \filldraw [black] (-1,0) circle (1pt) node[anchor=north]{\footnotesize{$-\varphi_1^\infty$}};

    \draw[lightgray,<-] (-1.2,1) -- (-1.3,1);    
    \draw[lightgray,<-] (1.2,1) -- (1.3,1);    

    \draw[<-] [domain=0.5:0.4,variable=\t,smooth] plot ({cos(\t r)},{sin(\t r)}); 
    \draw[<-] [domain=2.59:2.69,variable=\t,smooth] plot ({cos(\t r)},{sin(\t r)}); 
 
    \end{tikzpicture}
    \hspace{1.5cm}
    \begin{tikzpicture}[scale=1]
    \node (2) at (0, 2.5) {\boxed{\text{\footnotesize{$\lambda\in (\sigma^*_1,\sigma_1)$}}}};    
    
    \draw[->] (-1.5,0) -- (1.5,0) node[right] {\footnotesize{$U$}};
    \draw[->] (0,-0.1) -- (0,1.5) node[above] {\footnotesize{$z$}};

    \draw [domain=0:3.14,variable=\t,smooth] plot ({cos(\t r)},{sin(\t r)}); 
    
    \draw[lightgray,-] (-1.5,1) -- (1.5,1);

    \filldraw [blue] (0,1) circle (1pt) node[anchor=south west]{\footnotesize{$0$}};
    \filldraw [black] (1,0) circle (1pt) node[anchor=north]{\footnotesize{$\varphi_1^\infty$}};
    \filldraw [black] (-1,0) circle (1pt) node[anchor=north]{\footnotesize{$-\varphi_1^\infty$}};

    \filldraw [lightgray] (1,1) circle (1pt) node[anchor=south]{\footnotesize{$u^1$}};
    \filldraw [lightgray] (-1,1) circle (1pt) node[anchor=south]{\footnotesize{$-u^1$}};

    \filldraw [blue] (0.707,0.707) circle (1pt) node[anchor=west]{\footnotesize{$\mathcal{P}(u^1)$}};  
    \filldraw [blue] (-0.707,0.707) circle (1pt) node[anchor=east]{\footnotesize{$\mathcal{P}(-u^1)$}};

    \draw[lightgray,<-] (-1.2,1) -- (-1.3,1);    
    \draw[lightgray,<-] (1.2,1) -- (1.3,1);    

    \draw[lightgray,->] (-0.7,1) -- (-0.8,1);    
    \draw[lightgray,->] (0.7,1) -- (0.8,1);    

    \draw[<-] [domain=0.5:0.4,variable=\t,smooth] plot ({cos(\t r)},{sin(\t r)}); 
    \draw[<-] [domain=2.59:2.69,variable=\t,smooth] plot ({cos(\t r)},{sin(\t r)}); 

    \draw [blue,domain=0.8:2.4,variable=\t,smooth] plot ({cos(\t r)},{sin(\t r)}); 

    \draw[blue,->] [domain=1.2:1.1,variable=\t,smooth] plot ({cos(\t r)},{sin(\t r)}); 
    \draw[blue,->] [domain=1.85:1.95,variable=\t,smooth] plot ({cos(\t r)},{sin(\t r)}); 
 
    \end{tikzpicture}
    \hspace{1.5cm}
    \begin{tikzpicture}[scale=1]
    \node (1) at (0, 2.5) {\boxed{\text{\footnotesize{$\lambda\in (\sigma_1,\sigma^*_2)$}}}};    

    \draw[->] (-1.5,0) -- (1.5,0) node[right] {\footnotesize{$U$}};
    \draw[->] (0,-0.1) -- (0,1.5) node[above] {\footnotesize{$z$}};

    \draw [blue,domain=0:3.14,variable=\t,smooth] plot ({cos(\t r)},{sin(\t r)}); 
    
    \draw[lightgray,-] (-1.5,1) -- (1.5,1);

    \filldraw [blue] (0,1) circle (1pt) node[anchor=south west]{\footnotesize{$0$}};
    \filldraw [blue] (1,0) circle (1pt) node[anchor=north]{\footnotesize{$\varphi_1^\infty$}};
    \filldraw [blue] (-1,0) circle (1pt) node[anchor=north]{\footnotesize{$-\varphi_1^\infty$}};

    \draw[lightgray,->] (-1.2,1) -- (-1.3,1);    
    \draw[lightgray,->] (1.2,1) -- (1.3,1);    

    \draw[blue,->] [domain=0.5:0.4,variable=\t,smooth] plot ({cos(\t r)},{sin(\t r)}); 
    \draw[blue,->] [domain=2.59:2.69,variable=\t,smooth] plot ({cos(\t r)},{sin(\t r)}); 
 
    \end{tikzpicture}
    \begin{tikzpicture}[scale=1.5]
    \node (2) at (0, 2.25) {\boxed{\text{\footnotesize{$\lambda\in (\sigma^*_2,\sigma_2)$}}}};    

    \fill [blue!15, domain=0:3.14, variable=\t]
    (0, 0)
    -- plot ({cos(\t r)},{sin(\t r)})
    -- (0, 0)
    -- cycle;

    \filldraw[white,domain=0:3.14,variable=\t,smooth] plot ({cos(\t r)},{0.6*(sin(-7+\t r))^3}); 
    \draw[blue,domain=0:3.14,variable=\t,smooth] plot ({cos(\t r)},{0.6*(sin(-7+\t r))^3}); 
    \draw[blue,domain=2.1:2.13,variable=\t,smooth,->] plot ({cos(\t r)},{0.6*(sin(-7+\t r))^3}); 
    \draw[blue,domain=1.27:1.3,variable=\t,smooth,-<] plot ({cos(\t r)},{0.6*(sin(-7+\t r))^3}); 

    \draw[blue,dashed,domain=1.57:3.14,variable=\t,smooth] plot ({cos(\t r)},{0.8*(sin(4+\t r))^2}); 
    \draw[blue,dashed,domain=0:1.57,variable=\t,smooth] plot ({cos(\t r)},{0.8*(sin(4+\t r))^3}); 
    \draw[blue,domain=2.28:2.29,variable=\t,smooth,->] plot ({cos(\t r)},{0.8*(sin(4+\t r))^2}); 
    \draw[blue,domain=1:1.01,variable=\t,smooth,-<] plot ({cos(\t r)},{0.8*(sin(4+\t r))^3}); 
    

    \draw[-,dashed] (-1,0) -- (1,0);
    \draw[-] (-1,0) -- (-1.5,0);
    \draw[->] (1,0) -- (1.5,0);
    \draw[-,dashed] (0,-0.1) -- (0,1);
    \draw[->] (0,1) -- (0,1.5) node[above] {\footnotesize{$z$}};
    \draw[-,dashed] (-0.2,-0.2) -- (0.4,0.4);    
    \draw[->] (-0.2,-0.2) -- (-0.5,-0.5);

    \draw [blue,domain=0:3.14,variable=\t,smooth] plot ({cos(\t r)},{sin(\t r)}); 

    \draw[very thick] (-1,0) arc (180:360:1cm and 0.2cm);
    \draw[very thick,dashed] (-1,0) arc (180:0:1cm and 0.2cm);

    \draw[blue,rotate=-90,->] (-1,0) arc (180:230:1cm and 0.2cm);
    \draw[rotate=-90,-<] (-1,0) arc (180:240:1cm and 0.2cm);
    \draw[blue,rotate=-90,dashed,->] (-1,0) arc (180:147:1cm and 0.2cm);
    \draw[rotate=-90,dashed,-<] (-1,0) arc (180:136:1cm and 0.2cm);

    \draw[rotate=-90] (-1,0) arc (180:281:1cm and 0.2cm);
    \draw[blue,rotate=-90] (-1,0) arc (180:230:1cm and 0.2cm);
    \draw[rotate=-90,dashed] (-1,0) arc (180:101:1cm and 0.2cm);
    \draw[blue,rotate=-90,dashed] (-1,0) arc (180:140:1cm and 0.2cm);

    \draw[lightgray,-] (-0.75,1.25) -- (1.5,1.25);
    \draw[lightgray,-] (-1.5,0.5) -- (0.75,0.5);
    \draw[lightgray] (1.5,1.25) -- (0.75,0.5); 
    \draw[lightgray] (-0.75,1.25) -- (-1.5,0.5); 

    \draw[lightgray,<->] (-1,1) -- (1.25,1);
    \draw[lightgray,<->] (0.15,1.25) -- (-0.4,0.5); 

    \filldraw [lightgray] (-0.17,0.8) circle (1pt);
    \filldraw [lightgray] (0.07,1.15) circle (1pt);

    \filldraw [blue] (-0.16,0.6) circle (1pt); 
    \filldraw [blue] (0.13,0.8) circle (1pt); 




    \draw[blue,->>] [domain=0.5:0.4,variable=\t,smooth] plot ({cos(\t r)},{sin(\t r)}); 
    \draw[blue,->>] [domain=2.59:2.69,variable=\t,smooth] plot ({cos(\t r)},{sin(\t r)}); 

    \filldraw [blue] (0,1) circle (1pt) node[anchor=south west]{$0$};
    \filldraw [blue] (1,0) circle (1pt) node[anchor=north west]{\footnotesize{$\varphi_1^\infty$}};
    \filldraw [blue] (-1,0) circle (1pt) node[anchor=north east]{\footnotesize{$-\varphi_1^\infty$}};
    \filldraw [black] (-0.2,-0.2) circle (1pt) node[anchor=north west]{\footnotesize{$\varphi_2^\infty$}};
    \filldraw [black] (0.2,0.2) circle (1pt) node[anchor=south west]{\footnotesize{$-\varphi_2^\infty$}};    
 

    \end{tikzpicture}
    \hspace{1cm}
    \begin{tikzpicture}[scale=1.5]
    \node (2) at (0, 2.25) {\boxed{\text{\footnotesize{$\lambda>\sigma_2$}}}};    

    \fill [blue!15, domain=0:3.14, variable=\t]
    (0, 0)
    -- plot ({cos(\t r)},{sin(\t r)})
    -- (0, 0)
    -- cycle;
    
    \filldraw[blue!15, thick] (-1,0) arc (180:360:1cm and 0.175cm);

    \draw[-,dashed] (-1,0) -- (1,0);
    \draw[-] (-1,0) -- (-1.5,0);
    \draw[->] (1,0) -- (1.5,0);
    \draw[-,dashed] (0,-0.1) -- (0,1);
    \draw[->] (0,1) -- (0,1.5) node[above] {\footnotesize{$z$}};
    \draw[-,dashed] (-0.2,-0.2) -- (0.4,0.4);
    \draw[->] (-0.2,-0.2) -- (-0.5,-0.5);

    \draw [blue,domain=0:3.14,variable=\t,smooth] plot ({cos(\t r)},{sin(\t r)}); 

    \draw[blue,very thick] (-1,0) arc (180:360:1cm and 0.2cm);
    \draw[blue,very thick,dashed] (-1,0) arc (180:0:1cm and 0.2cm);
    \draw[blue,rotate=-90] (-1,0) arc (180:281:1cm and 0.2cm);
    \draw[blue,rotate=-90,dashed] (-1,0) arc (180:101:1cm and 0.2cm);

    \draw[blue,rotate=-90] (-1,0) arc (180:272:0.8cm and 0.105cm);
    \draw[blue] (1,0) arc (0:-80:0.8cm and 0.09cm);
    \draw [blue, domain=1.9:11,variable=\t,smooth,->] plot ({0.05*\t-0.2},{1/(\t)^2-0.1});

    \draw[blue,rotate=-90] (-1,0) arc (180:250:1cm and 0.295cm);
    \draw[blue] (-1,0) arc (180:256:0.4cm and 0.09cm);
    \draw [blue, domain=1.5:10,variable=\t,smooth,->] plot ({-0.05*\t-0.2},{1/(\t)^2-0.1});

    \draw[blue,->>] [domain=0.5:0.4,variable=\t,smooth] plot ({cos(\t r)},{sin(\t r)}); 
    \draw[blue,->>] [domain=2.59:2.69,variable=\t,smooth] plot ({cos(\t r)},{sin(\t r)}); 

    \draw[blue,rotate=-90,->] (-1,0) arc (180:250:1cm and 0.2cm);
    \draw[blue,rotate=-90,dashed,->] (-1,0) arc (180:115:1cm and 0.2cm);

    \draw[blue,very thick,-<] (-1,0) arc (180:240:1cm and 0.2cm);
    \draw[blue,very thick,dashed,-<] (-1,0) arc (180:115:1cm and 0.2cm);    

    \draw[blue,very thick,->] (-1,0) arc (180:280:1cm and 0.2cm);
    \draw[blue,very thick,dashed,->] (-1,0) arc (180:44:1cm and 0.2cm);    

    \draw[lightgray,-] (-0.75,1.25) -- (1.5,1.25);
    \draw[lightgray,-] (-1.5,0.5) -- (0.75,0.5);
    \draw[lightgray] (1.5,1.25) -- (0.75,0.5); 
    \draw[lightgray] (-0.75,1.25) -- (-1.5,0.5); 

    \draw[lightgray,<<->>] (-1,1) -- (1.25,1);
    \draw[lightgray,<->] (0.15,1.25) -- (-0.4,0.5); 
    
    \filldraw [blue] (-0.02,1) circle (1pt) node[anchor=south west]{$0$};
    \filldraw [blue] (1,0) circle (1pt) node[anchor=north west]{\footnotesize{$\varphi_1^\infty$}};
    \filldraw [blue] (-1,0) circle (1pt) node[anchor=north east]{\footnotesize{$-\varphi_1^\infty$}};
    \filldraw [blue] (-0.2,-0.2) circle (1pt) node[anchor=north west]{\footnotesize{$\varphi_2^\infty$}};
    \filldraw [blue] (0.2,0.2) circle (1pt) node[anchor=south west]{\footnotesize{$-\varphi_2^\infty$}};

    \end{tikzpicture}
\caption{The global dynamics of the compactified semiflow as the parameter $\lambda$ increases and passes through the values $\sigma^*_1,\sigma_1,\sigma^*_2,\sigma_2$. As $\lambda$ crosses $\sigma^*_1$, the trivial equilibrium $u_*\equiv 0$ undergoes a pitchfork bifurcation, yielding the bounded and stable equilibria $\pm u^1(x)$, and their respective heteroclinics. As $\lambda\to \sigma_1$, such equilibria $\pm u^1(x)$ become arbitrarily large and its compactification, $\pm \mathcal{P}(u^1(x))$, converge to the Steklov eigenfunctions $\pm \Phi_1$. For $\lambda>\sigma_1$, there are stable equilibria at infinity, given by $\pm \varphi^\infty_1(x)$, that attract infinite-time blow-up solutions. Note that $\pm \varphi^\infty_1(x)=\pm \Phi_1(x)$ when $\lambda=\sigma_1$. Similarly as $\lambda$ crosses $\sigma^*_2$ and $\sigma_2$. The compactification of the attractor is depicted in blue.
}\label{Fig:P1d}
\end{figure}

\section{Discussion}\label{sec:disc}

We discuss on our approach of bifurcation from infinity and provide open problems.

Previously, in \cite{Rabinowitz,ARRB07,ARRB09}, it was known that there are unbounded bifurcation curves of equibria from infinity that converges, after rescaling, to the Steklov eigenfunctions as the parameter approaches a Steklov eigenvalue, as in~\eqref{EQconvSteklov}.
In our compactification approach, we show that such Steklov eigenfunctions are equilibria of the induced semiflow at infinity. 
We have also discussed the asymptotic behavior of solutions in Lemma~\ref{lem:nondiss}, which proves that there are infinite-time blow-up solutions that converge, after rescaling, to certain eigenfunctions (which are equilibria at infinity).
Therefore, the Poincaré projection justifies the interpretation of such infinite-time blow-up solutions as \emph{heteroclinics to infinity}. 
The projection also allows us to identify the induced behavior at infinity and compute the compactification of the unbounded attractor.

In particular, for $\lambda>\sigma_2$, note that there are sufficiently small initial data that shadow the heteroclinic from $0$ towards $\varphi_2^\infty$ for some time, but eventually converge to $\varphi_1^\infty$, see Figure~\ref{Fig:P1d}.
Such heteroclinics may cause a numerical approximation to break down for sufficiently large times, leading to the false conclusion that any sufficiently small solution may approximate $\varphi^\infty_2$, after rescaling. However, this is only a transient behavior, as rescaled solutions generically converge to $\varphi_1$. Therefore, our compactification approach may provide a more accurate and reliable approximation scheme for large solutions.

For simplicity, we have considered a scalar PDE in one spatial dimension with boundary nonlinearities $g(u)$ which are Lipschitz, odd, monotone, bounded, $g(0)=0$ and $g'(0)=1$; for example $g(u)=\arctan(u)$. 
However, if any of these properties is violated, different phenomena may occur. In some cases, our result can still be replicated, but with a slightly different outcome. 
For example, $g(0)=0$ guarantees that $u_*\equiv 0$ is an equilibrium of \eqref{intro:PDE}, but one can consider a different bounded equilibrium that undergoes pitchfork bifurcations and construct the resulting attractor. 
Alternatively, if $g'(0)=-1$, such as in the example $g(u)=-\arctan (u)$, the subcritical/supercritical character of the bifurcations at $u_*\equiv 0$ and from infinity change roles, see
Figure~\ref{Fig:Bif}; with a different dynamical structure of the attractor for each parameter value.
In general, for non-odd or non-monotone nonlinearities, additional bounded equilibria may arise, and thus one has to construct all the bounded heteroclinic connections, akin to \cite{LappicyFully} and references therein, and discuss the connection problem for the infinite-time blow-up solutions, as in~\cite{BenGal10}.

For example, we mention a few possible future explorations in this direction:
\begin{enumerate}
    \item Nonlinearities $g(\lambda,u)$ which are not necessarily monotone in $\lambda$, may not possess monotone bifurcating curves from infinity, 
    see \cite[Theorem 4.1 and 4.2]{ARRB09}. 
    \item If $g(u)\neq 0$ for any $u\in\mathbb{R}$, such as sigmoid function $g(u)=1/(1+\mathrm{e}^{-u})$, then the bifurcating branch from infinity may have to meet another bifurcation from infinity. 
    \item The case of separated nonlinear boundary conditions, such that each boundary condition has a nonlinearity given by $g_0(u),g_1(u)$ at $x=0,1$.
    \item In case that $g(x,u)$ has different linearizations at $u_*\equiv 0$ and `at infinity', in constrast to~\eqref{EVprob} and \eqref{EVprobinfty}. 
    \item \textcolor{black}{Parabolic PDEs with nonlinear boundary conditions in higher spatial dimensions, which amounts to a nonlinear Steklov problem on the unit circle; 
    see \cite{Cushing2,Cushing3}. 
    }
\end{enumerate}
%
%
Next, one may consider the interaction of the current nonlinear Robin boundary conditions with extra nonlinear terms in the interior of the domain. 
The role of different reaction terms, and an appropriate compactification, is provided in \cite{LP}; specially in the role of a bifurcation at infinity, see \cite{IchSaka24}. A natural next step is to analyze the interplay between the dynamics induced by the interior reaction and the boundary nonlinearities.



\textbf{Acknowledgments.} JMA was partially supported by grants PID2022-137074NB-I00 and CEX2023-001347-S “Severo Ochoa Programme for Centres of Excellence in R\&D” MCIN/AEI/10.13039/501100011033, the two of them from Ministerio de Ciencia e Innovación, Spain. Also by “Grupo de Investigaci\'on 920894 - CADEDIF”, UCM, Spain. JF was funded by CNPq 406460/2023-0. PL was supported by Marie Skłodowska-Curie Actions, UNA4CAREER H2020 Cofund, 847635, 
with the project DYNCOSMOS.

\textbf{Competing interest and data availability.}
The authors have no conflict of interest to declare. Moreover, the authors confirm that the data supporting the findings of this study are available within the article.


\medskip

\begin{thebibliography}{10}
\small{

\bibitem{ACRB99}
J.M.~Arrieta, A.N.~Carvalho and A.~Rodriguez-Bernal.
Parabolic problems with nonlinear boundary conditions and critical nonlinearities.
\emph{J. Diff. Eq.} \textbf{156}, 376--400, (1999).

\bibitem{AC00}
J.M.~Arrieta and A.N.~Carvalho.
Abstract parabolic problems with critical nonlinearities and applications to Navier-Stokes and heat equations.
\emph{Trans. A. M. S.} \textbf{352}, 285--310, (2000).

\bibitem{ACRB00}
J.M.~Arrieta, A.N.~Carvalho and A.~Rodriguez-Bernal.
Attractors of parabolic problems with nonlinear boundary conditions. Uniform bounds.
\emph{Comm. P. D. E.} \textbf{25}, 1--37, (2000).

\bibitem{ARRB07}
J.M.~Arrieta, R.~Pardo and A.~Rodriguez-Bernal.
Bifurcation and stability of equilibria with asymptotically linear boundary conditions at infinity.
\emph{Proc. Royal Soc. Edinburgh} \textbf{137A}, 225--252, (2007).

\bibitem{ARRB09}
J.M.~Arrieta, R.~Pardo and A.~Rodriguez-Bernal.
Equilibria and global dynamics of a problem with bifurcation from infinity.
\emph{J. Diff. Eq.} \textbf{246}, 2055--2080, (2009).

\bibitem{ARRB10}
J.M.~Arrieta, R.~Pardo and A.~Rodriguez-Bernal.
Infinite Resonant Solutions and Turning Points in a Problem with Unbounded Bifurcation.
\emph{Int. J. Bif. Chaos} \textbf{20}, 2885--2896, (2010).

\bibitem{BenGal10}
N.~Ben-Gal.
Grow-Up Solutions and Heteroclinics to Infinity for Scalar Parabolic PDEs.
\emph{Ph.D. Thesis, Division of Applied Mathematics, Brown University}, (2010).

\bibitem{FernandesBortolan}
M.~Bortolan and J.~Fernandes.
Sufficient conditions for the existence and uniqueness of maximal attractors for autonomous and nonautonomous dynamical systems. 
\emph{J. Dyn. Diff. Eq.} \textbf{1}, 1--30, (2022).

\bibitem{BCP}
S. Bruschi, A. N. Carvalho, J. Pimentel.
Limiting grow-up behavior for a one-parameter family of dissipative PDEs. 
\emph{Indiana Univ. Math. J.} \textbf{69}, 657--683, (2020)

\bibitem{CP12}
A.~Castro and R.~Pardo.
Resonant Solutions and Turning Points in an Elliptic Problem with Oscillatory Boundary Conditions.
\emph{Pacific J. Math.} \textbf{257}, (2012).

\bibitem{CP17}
A.~Castro and R.~Pardo.
Infinitely Many Stability Switches in a Problem with Sublinear Oscillatory Boundary Conditions.
\emph{J. Dyn. Diff. Eq.} \textbf{29}, 485--499, (2017).

{\color{black}
\bibitem{Cushing}
J.M.~Cushing.
Some existence theorems for nonlinear eigenvalue problems associated with elliptic equations.
\emph{Arch. Rat. Mec. Anal.,} \textbf{42}, 63--76, (1971).

\bibitem{Cushing2}
J.M.~Cushing.
Nonlinear Steklov problems on the unit circle.
\emph{J. Math. Anal. Appl.,} \textbf{38}, 766--783, (1972).

\bibitem{Cushing3}
J.M.~Cushing.
Nonlinear Steklov problems on the unit circle II (and a Hydrodynamical Application).
\emph{J. Math. Anal. Appl.,} \textbf{39}, 267--378, (1973).

\bibitem{Dancer}
E.N.~Dancer.
A note on bifurcation from infinity. 
\emph{Quaterly J. Math.} \textbf{25}, 81--84, (1974).
}

\bibitem{RochaFiedler}
B.~Fiedler and C.~Rocha.
Heteroclinic orbits of semilinear parabolic equations. 
\emph{J. Diff. Eq.} \textbf{125}, 239--281, (1996).

\bibitem{Hell09}
J.~Hell.
Conley Index at Infinity. 
\emph{Topol. Methods Nonlin. Anal.} \textbf{42}, 137--167, (2013).

\bibitem{Henry81}
D.~Henry.
\emph {Geometric Theory of Semilinear Parabolic Equations}.
Springer-Verlag New York, (1981).

\bibitem{IchSaka24}
Y.~Ichida and T.O.~Sakamoto.
Geometric Approach to the Bifurcation at Infinity: A Case Study.
\emph{Qual. Theor. Dyn. Sys.} \textbf{23}, 109, (2024).

\bibitem{LappicyFully}
P.~Lappicy.
Sturm attractors for fully nonlinear parabolic equations in one spatial dimension.
\emph{Rev. Mat. Complutense} \textbf{36}, 725--747, (2023).

\bibitem{LP}
P.~Lappicy and J.~Fernandes.
Unbounded Sturm attractors for quasilinear parabolic equations.
\emph{Proc. Edinburgh Math. Soc.} \textbf{67} 542 -- 565, (2024).

\bibitem{LBeatriz}
P.~Lappicy and E.~Beatriz.
An energy formula for fully nonlinear degenerate parabolic equations in one spatial dimension.
\emph{Math. Ann.} \textbf{389} 4125 -- 4147, (2024).

\bibitem{P}
J.~Pimentel.
Unbounded Sturm global attractors for semilinear parabolic equations on the circle. 
\emph{SIAM Journal on Mathematical Analysis} \textbf{48} 3860 -- 3882, (2016).

\bibitem{RochaPimentel15}
J.~Pimentel and C.~Rocha.
Noncompact global attractors for scalar reaction–diffusion equations.
\emph{Sao Paulo J. Math.} \textbf{9}, 299--310, (2015).

\bibitem{RochaPimentel16}
J.~Pimentel and C.~Rocha.
A permutation related to non-compact global attractors for slowly non-dissipative systems.
\emph{J. Dyn. Diff. Eq.} \textbf{28}, 1--15, (2016).

\bibitem{Rabinowitz}
P.~Rabinowitz.
On bifurcation from infinity.
\emph{J. Diff. Eq.} \textbf{14}, 462--475, (1973).

{\color{black}
\bibitem{Stekloff}
W.~Stekloff.
Sur les problèmes fondamentaux de la physique mathématique.
\emph{Annales scientifiques de l’É.N.S., 3e série,} \textbf{19}, 455--490, (1902).

\bibitem{Stuart}
C.A.~Stuart.
Solutions of large norms for nonlinear Sturm-Liouville problems.
\emph{Q. J. Math.} \textbf{2}, 129--139, (1973).

\bibitem{StuartToland}
C.A.~Stuart and J.F.~Toland.
A Global Result Applicable to Nonlinear Steklov Problems.
\emph{J. Diff. Eq.} \textbf{15}, 247--368, (1974).

\bibitem{Toland}
J.F.~Toland.
Asymptotic linearity and nonlinear eigenvalue problems.
\emph{Q. J. Math.} \textbf{2}, 241--250, (1973).

\bibitem{Toland2}
J.F.~Toland.
Asymptotic linearity and nonlinear eigenvalue problems II.
\emph{Proc. R. IR. Acad. Sect. A} \textbf{77}, 1--12, (1977).
} 

}
\end{thebibliography}
\end{document}